\documentclass[journal]{IEEEtran}
\usepackage{blindtext}
\usepackage{graphicx}
\usepackage{balance}
\usepackage[utf8]{inputenc}
\usepackage{mathtools}
\usepackage{amsmath,mathtools,amssymb,amsfonts,amsthm}

\usepackage{tikz, subcaption, cite}
\usepackage{tabu}
\usepackage{}

\newtheorem{lemma}{Lemma}
\newtheorem{remark}{Remark}

\usepackage{epsfig} 
\usepackage{psfrag}
\usepackage{mathrsfs}
\usepackage{nomencl}   
\makenomenclature
\usepackage{multirow}
\usepackage{multicol}
\usepackage{array,booktabs,arydshln,xcolor}
\usepackage{soul}

\usepackage{mathtools}
\usepackage{nccmath}
\usepackage{caption}
\usepackage{nccmath}
\usepackage{comment}
\usepackage{makecell}
\graphicspath{{Figs/}}

\usepackage{siunitx}
\usepackage{multirow}

\usepackage{tikz}
\usetikzlibrary{fit,positioning,automata,backgrounds}
\usetikzlibrary{calc} 

\newlength{\figW}   
\usepackage{stackengine}

\DeclareMathOperator*{\argmax}{arg\,max}  
\DeclareMathOperator*{\argmin}{arg\,min}  

\usepackage{amsmath}
\usepackage{amssymb}
\usepackage{float}
\usepackage{dblfloatfix} 
\usepackage{amsthm}
\usepackage{graphicx}
\usepackage{epstopdf}
\usepackage{color}
\usepackage{verbatim}
\usepackage{tikz,times}
\usepackage{xcolor}

\usepackage{hyperref}
\hypersetup{ colorlinks=true, linkcolor=blue, filecolor=magenta, urlcolor=cyan,}

\usepackage{enumerate}

\usepackage[ruled,vlined]{algorithm2e}

\newcommand{\matrixstyle}[1]{{#1}}

\newcommand{\expect}[1]{\mathbb{E}_{#1}}
\newcommand{\entropy}[1]{\mathbb{H}\left[#1\right]}

\newcommand{\trace}{\mathrm{tr}}

\newtheorem{assumption}{\textbf{Assumption}}

\renewcommand{\baselinestretch}{0.93}

\makeatletter
\renewcommand*\env@matrix[1][*\c@MaxMatrixCols c]{%
  \hskip -\arraycolsep
  \let\@ifnextchar\new@ifnextchar
  \array{#1}}
\makeatother

\ifCLASSINFOpdf
\else
\fi

\begin{document}

\title{\huge Gaussian Process Dual MPC using Active Inference: \\ An Autonomous Vehicle Usecase}

\author{Mohammad Mahmoudi Filabadi, Guillaume Crevecoeur and Tom Lefebvre
\thanks{M. M. Filabadi, G. Crevecoeur and T. Lefebvre are with the Department of Electromechanical, Systems and Metal Engineering, Ghent University, B-9052 Ghent, Belgium, e-mail: 
{\tt\small \{mohammad.mahmoudifilabadi, tom.lefebvre, guillaume.crevecoeur\}@ugent.be}.}
\thanks{M. M. Filabadi, G. Crevecoeur and T. Lefebvre are member of core lab MIRO, Flanders Make, Belgium.}
}
\maketitle
\begin{abstract}

Designing controllers under uncertainty requires balancing the need to explore system dynamics with the requirement to maintain reliable control performance. Dual control addresses this challenge by selecting actions that both regulate the system and actively gather informative data. 
This paper investigates the use of the Active Inference framework, grounded in the Free Energy Principle, for developing a dual model-predictive controller (MPC).
To identify and quantify uncertainty, we introduce an online sparse semi-parametric Gaussian Process model that combines the flexibility of nonparametric with the efficiency of parametric learning for real-time updates.
By applying the expected free energy functional to this adaptive probabilistic model, we derive an MPC objective that incorporates an information-theoretic term, which captures uncertainty arising from both the learned model and measurement noise. 
This formulation leads to a stochastic optimal control problem for dual controller design, which is solved using a novel dynamic-programming-based method. 
Simulation results on a vehicle use case demonstrate that the proposed algorithm enhances autonomous driving control performance across different settings and scenarios.

\end{abstract}

\section{Introduction} \label{sec:Intro}

Accurate control of dynamical systems often relies on precise knowledge of system dynamics, yet in many practical applications, these dynamics are only partially known or subject to uncertainty. 
Traditional model predictive control (MPC) frameworks typically assume that a reasonably accurate model is already available, leading to suboptimal or unstable performance when this assumption is violated. 
To address this limitation, dual control strategies have been developed in which the controller not only aims to meet performance goals but also actively probes the system to improve model accuracy over time  \cite{Mesbah2018, Heirung2017, Li2025, Alpcan2015}.
{
The fundamental dual control problem, formulated as a Partially Observable Markov Decision Process (POMDP), entails solving an SOC problem in the belief space.
In this formulation, the parametric uncertainties act as latent variables, and information about them must be inferred indirectly through state measurements. 
Closed-form solutions for this POMDP are generally intractable. Existing approaches address this challenge either by approximating the problem, known as implicit dual control (IDC), or by explicitly incorporating probing behavior into the control inputs, referred to as explicit dual control (EDC). EDC methods typically employ heuristic strategies, either by directly injecting system excitation or by reformulating the decision-making problem to account for model uncertainty over future control stages \cite{Mesbah2018}.
}
Recent works have extended the notion of EDC using information-theoretic formulations \cite{Alpcan2015}, where the controller maximizes information gain about the uncertain system while optimizing performance. 

This approach to dual control reflects the broader challenge of motivating artificial agents to acquire useful information, a problem that remains open across many disciplines.
A promising framework that has attracted increasing attention for addressing this issue is the Active Inference Framework (AIF), a paradigm originating in theoretical neuroscience \cite{parr2022active, smith2022step} that provides a unified account of perception, action, and learning by biological agents such as the human brain. 
AIF proposes that agents maintain an internal generative model, a probabilistic representation of how the world behaves and how their actions affect it. Using this model, they act to minimize expected free energy (EFE), a quantity that represents the anticipated uncertainty or surprise associated with future sensory information \cite{Millidge2021, Kouw2024}.  
By minimizing EFE, agents pursue actions that naturally balance epistemic (information-seeking) value---capturing how much a policy is expected to reduce uncertainty about latent states---and extrinsic (goal-seeking) value, which measures how well the policy’s expected observations align with the agent’s goals. 
{
In this way, AIF provides a rigorous foundation for information-driven EDCs despite their heuristic nature.
}

Several studies have already employed AIF in the fields of robotics and control engineering, including \cite{Pezzato2023, PioLopez2016, Meo2023, Pezzato2020, Meo2021, baioumy2021active, Kouw2024}. Particularly, \cite{PioLopez2016, Meo2023, Pezzato2020, Meo2021, baioumy2021active} employ a simplified formulation of free energy minimization that incorporates both the mean-field and Laplace approximations. These studies define a reference model to specify the desired behavior of the system as the generative model.
The recent study by \cite{Kouw2024} introduced a new approach by employing a parametric model learner as the generative model and applying exact EFE minimization. However, it still relies on the mean-field approximation and disregards statistical dependencies in the resulting optimization problem by collapsing the posterior predictive distributions to their most probable values. 

In this paper, we address the limitations of previous studies by extending EFE minimization to a more general framework for probabilistic model learning. Our model learning approach is based on an online sparse Gaussian Process (GP) formulation that accommodates both parametric and nonparametric model components. 
{
By minimizing the EFE, an identification objective, derived from the information gain of the GP model, is superposed  into the task related cost. 
As a result, the control actions are selected to strategically explore the system, gathering informative data that improves long-term control performance.}
In addition, we consider the statistical dependencies among random variables by formulating the problem as one of SOC. Finally, we present a probabilistic trajectory optimization method based on iterative local approximation of dynamic programming (DP) to solve the resulting SOC problem. 
The preliminary results of the proposed probabilistic trajectory optimization approach are reported in \cite{Filabadi2025}.

GP models offer a probabilistic, nonparametric approach for black-box system identification. They efficiently exploit available data and quantify uncertainty. GPs are widely applied in robotics and control, where modelling uncertainty is crucial for decision-making \cite{Pan2014, Deisenroth2015, Boedecker2014, Frigola2014, Filabadi2025, Hashimoto2025, He2024,  Chowdhary2015, Benciolini2024}. 
Despite the high learning efficiency of GP models, their computational cost scales cubically with the number of training points \cite{Liu2020_GPBigData}, making them problematic for online applications where data are acquired sequentially \cite{He2024, Chowdhary2015}.

This work proposes a semi-parametric extension of the online sparse GP \cite{Csató2002, He2024, Chowdhary2015}, combining the flexibility of non-parametric GPs with the interpretability of parametric models. Online sparse GPs maintain a sparse set of representative data points, reducing the computational and memory burden typically associated with standard GPs and enabling real-time updates as new data become available. We demonstrate that incorporating a parametric component into online sparse GP models enhances learning efficiency and reduces data requirements, with the parametric model capturing information from all observed data while the non-parametric component preserves a compact, informative subset of data.


The concept of generating informative data by maximizing information gain in GP model learning has also been widely employed for active dynamics exploration and efficient system identification \cite{Patwardhan2014, Memmel2023, Bayen2020, yu2021, Sukhija2023, Benciolini2024}. Furthermore, the simultaneous combination of active learning and control design has also been investigated in \cite{Alpcan2011, Le2021, Ladislav2014}. However, these works generally overlook the computational challenges associated with real-time implementation and do not account for adaptive learning mechanisms. In addition, the optimal dual control naturally emerges from solving an SOC problem when it is formulated particularly under GP dynamics. 
This aspect is often neglected in prior studies, which typically consider only the GP predictive mean and fail to incorporate the predictive covariance (uncertainty) into the decision-making process. 

This paper examines the exact SOC formulation derived from EFE minimization applied to the model identified via the aforementioned adaptive learning scheme. 
We optimize the resulting objective function using a receding-horizon strategy that simultaneously identifies system dynamics, quantifies uncertainty, and optimizes control actions through a probabilistic trajectory optimization. 
Finally, we validate our approach on a vehicle dynamics model across various tasks and settings, enabling the development of an intelligent motion planner for autonomous vehicles.

In the following sections, we first present the problem formulation in Sec. \ref{sec:Preliminaries}. The subsequent sections outline the proposed approach, with each section highlighting its respective contribution as follows:
\begin{itemize}
    \item
    Sec. \ref{sec:Inf4Model}: Introduces a novel semi-parametric extension of online sparse GP as the generative model for AIF.
    
    \item
    Sec. \ref{sec:ActInf}: Provides a rigorous foundation for EDC by applying EFE minimization to the proposed generative model, resulting in a stochastic dynamic optimization problem that balances information-seeking and goal-directed behavior. 
    
    \item
    Sec. \ref{sec:Traj_Opt}: Presents a probabilistic trajectory optimization method to solve the formulated stochastic dynamic optimization problem.
    
    \item
    Sec. \ref{sec:algorithm}: Integrates these components to propose an adaptive dual MPC algorithm.
    
    \item
    Sec. \ref{sec:Results}: Validates the proposed methods and algorithm through simulation studies, including applications to autonomous vehicle development.

    \item
    Finally, Sec. \ref{sec:Conclusion} provides conclusions and discusses potential directions for future work.
\end{itemize}

\section{Preliminaries} \label{sec:Preliminaries}

\subsection{Notation}
The notation $p(A)$ denotes the probability density function of $A$, while $p(A|B)$ represents the conditional probability of $A$ given $B$. 
The notation $\mathcal{N}(X;m,P)$ signifies the normal or Gaussian distribution of the random variable $X$ with mean $m$ and a positive-definite covariance matrix $P$.  
The operator \( \mathbb{E}[\cdot] \) is used to represent the expected value of a random variable. When this operator includes a subscript, it specifies the probability measure or random variable with respect to which the expectation is computed.
The entropy of the random variable $A$ is denoted as $\entropy{A}$ or $\entropy{p(A)}$ and is defined as $\expect{p(A)}[-\ln p(A)]$.
We use a white square, $\square$, to indicate symmetric entries within matrices. Finally, $|\mathcal{D}|$ denotes the cardinal number of the set $\mathcal{D}$. 


\subsection{Problem Formulation} \label{sec:Problem_Formulation}
We consider a class of nonlinear discrete-time stochastic systems that are governed by the following difference equation.
\begin{equation} \label{eq:general_dynamics_transition}
    x_{t+1} = f(x_t, u_t) + v_t, \quad t = 0,1,2,\dots,
\end{equation}
Here, $x_{t}\in \mathbb{R}^{n_x}$ and $u_{t}\in \mathbb{R}^{n_u}$ represent the state of the dynamic system with components $x_{i,t}$ for $i = 1, \dots, n_x$ and action vector at the discrete time $t$, respectively. The stochasticity of the system is modeled by the zero-mean Gaussian noise $v_{t} \in \mathbb{R}^{n_x}$ distributed by $\mathcal{N}\left( v_{t} ; 0, \Sigma_f \right)$ in which $\Sigma_f \in \mathbb{R}^{n_x \times n_x}$. 
We define the concatenated state-action vector at time $t$ as $\xi_t \coloneqq \begin{bmatrix}
    x_t^{\top} & u_t^{\top}
\end{bmatrix}^{\top} \in \mathbb{R}^{n_{\xi}}$ with $n_{\xi} \coloneqq n_x+n_u$. The initial state $x_0$ is assumed to be drawn from a known normal distribution $\mathcal{N}\left(x_{0} ; \mu_{x,0}, \Sigma_{xx,0} \right)$. We also assume that the states are fully observable.

In this work, we employ a stochastic MPC framework \cite{Christophersen2007}. 
Specifically, at each time step, $t$, we solve an optimization problem over a finite horizon $H$, subject to the uncertain dynamics (\ref{eq:general_dynamics_transition}), to determine the optimal action trajectory according to a given optimality criterion. Then, the first element of the optimal action trajectory will be applied to the system.  
In this formulation, the system dynamics (\ref{eq:general_dynamics_transition}) is expressed as a controlled Markov model with transition probability defined by
\begin{equation}
     p(x_{k+1}|\xi_k) = \mathcal{N}(x_{k+1};f(\xi_k), \Sigma_f),
\end{equation}
for $k=t,\dots,t+H-1$.
Consider the following set of random variables at time $t$ over a finite horizon $H$: 
\begin{equation}
    \begin{aligned}  
        U_t & \coloneqq \{u_t, u_{t+1}, \dots, u_{t+H-1} \}, \\  
        \bar{X}_t & \coloneqq \{x_t, x_{t+1}, \dots, x_{t+H-1} \}, \\
        \bar{Z}_t & \coloneqq \{\xi_t, \xi_{t+1}, \dots, \xi_{t+H-1} \},
    \end{aligned} 
\end{equation}
and $X_t \coloneqq \bar{X}_t \cup \{ x_{t+H} \}$, $Z_t \coloneqq \bar{Z}_t \cup \{ x_{t+H} \}$, such that $\bar{Z}_t = (U_t, \bar{X}_t)$ and $Z_t = (U_t, X_t)$.
Then, we describe the probabilistic dynamics over the horizon through the following factorization.
\begin{equation}\label{eq:factorization}
p(X_t|U_t) = p(x_t)\prod\nolimits_{k = t}^{t+H-1} p(x_{k+1}|\xi_k).
\end{equation}
Finally, we consider the following objective function at each time step:  
\begin{equation} \label{eq:SOC_objective}  
    J_{t}(U_t) \coloneqq  \mathbb{E}_{p(X_t|U_t)} \left[ C_t\left( Z_t \right) \right],  
\end{equation}  
where $C_t(Z_t)$ is a standard cumulative cost defined as
\begin{equation}\label{eq:cumulative_cost}
    C_t(Z_t) \coloneqq \sum\nolimits_{k=t}^{t+H} c_{k}(\xi_k),
\end{equation}
Here, $c_k:\mathbb{R}^{n_\xi}\mapsto \mathbb{R}^+$ denotes a non-negative stage (or immediate) cost function, while the terminal cost $c_{t+H}(\xi_{t+H})$ depends only on the terminal state $x_{t+H}$, i.e., $c_{t+H}(\xi_{t+H}) \equiv c_{t+H}(x_{t+H})$. 
We aim to find the optimal sequence $U^*_t$ such that 
\begin{equation}\label{eq:Optimization_Problem}
    U^*_t = \argmin_{U_t} J(U_t),
\end{equation}
in the presence of uncertain dynamics (\ref{eq:general_dynamics_transition}). 

The unknown state-space function $f:\mathbb{R}^{n_x}\times\mathbb{R}^{n_u} \mapsto \mathbb{R}^{n_x}$ is assumed to be a deterministic Lipschitz continuous nonlinear mapping and its components corresponding to each subsystems of the dynamic system are given by $f_i:\mathbb{R}^{n_x}\times\mathbb{R}^{n_u} \mapsto \mathbb{R}$ for $i = 1, \dots, n_x$. 
Since $f$ is unknown, we can also consider its value as a random variable.
To that end, we have to provide the estimation of the function $f$ over the planning horizon. This estimation will be represented by the following set of random variables
\begin{equation}
    \hat{F}_t  \coloneqq \{\hat{f}_{t+1}, \hat{f}_{t+2}, \dots ,\hat{f}_{t+H} \}.
\end{equation}

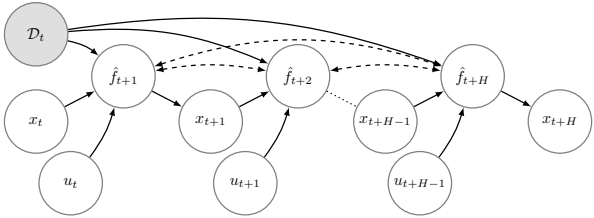
\begin{figure}[!h] 
	\centering 
	\resizebox{0.9\columnwidth}{!}{
    \begin{tikzpicture}
        \tikzstyle{main}=[circle, minimum size = 14mm, thick, draw =black!50, node distance = 5mm]
        \tikzstyle{connect} = [-latex,thick]
        \tikzstyle{floating} = [thick,dotted]
        \tikzstyle{curve} = [-latex,thick,bend right=10]
        \tikzstyle{curve2} = [-latex,thick,bend left=14]

        \tikzstyle{curve3} = [latex-latex, thick, dashed, bend left=10]
        \tikzstyle{curve4} = [latex-latex, thick, dashed, bend left=18]

        \node[main] (x_t) {$x_t$};
        \node[main] (u_t) [below right=of x_t, xshift=-0.6cm] {$u_t$};
        \node[main] (f_t1) [right=of x_t, yshift=10mm] {$\hat{f}_{t+1}$};
        \node[main] (x_t1) [right=of f_t1, yshift=-10mm] {$x_{t+1}$};
        \node[main] (u_t1) [below right=of x_t1, xshift=-0.6cm] {$u_{t+1}$};
        \node[main] (f_t2) [right=of x_t1, yshift=10mm] {$\hat{f}_{t+2}$};
        \node[main] (x_tH1) [right=of f_t2, yshift=-10mm,label=center:$x_{t+H-1}$] {};
        \node[main] (u_tH1) [below right=of x_tH1, xshift=-0.6cm, label=center:$u_{t+H-1}$] {};
        \node[main] (f_tH1) [right=of x_tH1, yshift=10mm,label=center:$\hat{f}_{t+H}$] {};
        \node[main] (x_tH) [right=of f_tH1, yshift=-10mm,label=center:$x_{t+H}$] {};
        
        
        \node[main,fill=gray!25] (D) [above = of x_t] {$\mathcal{D}_t$};
        
        \path (u_t) edge [curve] (f_t1);
        \path (u_t1) edge [curve] (f_t2);
        \path (u_tH1) edge [curve] (f_tH1);
    
        \path (x_t) edge [connect] (f_t1);
        \path (x_t1) edge [connect] (f_t2);
        \path (x_tH1) edge [connect] (f_tH1);
    
        \path (f_t1) edge [connect]  (x_t1);
        \path (f_t2) edge [floating] (x_tH1);
        \path (f_tH1) edge [connect] (x_tH);
    
        \path (D) edge [curve2] (f_t1);
        \path (D) edge [curve2] (f_t2);
        \path (D) edge [curve2] (f_tH1);

        \path (f_t1) edge [curve3] (f_t2);
        \path (f_t2) edge [curve3] (f_tH1);
        \path (f_t1) edge [curve4] (f_tH1);

    \end{tikzpicture}
	}
	\caption{\small{The probabilistic model used to represent the dependencies between variables over a finite horizon.}}\label{fig:PGM}
\end{figure}

The implication of the fact that $f$ is unknown on the MPC problem formulation above will be discussed and detailed in Section \ref{sec:ActInf}.
Furthermore, to quantify uncertainty and estimate the system's dynamics at time $t$, we use real-time data collected up to time $t$ in the following dataset:
\begin{equation} \label{eq:data_set_D}
    \mathcal{D}_t \coloneqq \{ (\xi_0,x_1), (\xi_1,x_2), \dots, (\xi_{t-1}, x_t)  \}. 
\end{equation}
This set can be split into the following datasets.
\begin{equation} \label{eq:data_set_D_split}
    \begin{aligned}
        \mathcal{D}^{\xi}_t & \coloneqq \{\xi_0, \xi_1, \dots, \xi_{t-1}\}, \\
        \mathcal{D}^{x_i}_t & \coloneqq \{x_{i,1}, x_{i,2} \dots, x_{i,t}\}, \\
        \mathcal{D}^{i}_t & \coloneqq \{(\xi_0,x_{i,1}), (\xi_1,x_{i,2}), \dots, (\xi_{t-1}, x_{i,t})\}, \\
    \end{aligned}
\end{equation}
where $\mathcal{D}^{\xi}_t$ denotes the set of the feature vectors, and $\mathcal{D}^{x_i}_t$ represents the labels for learning the subsystem $f_i$. Both are collected in $\mathcal{D}^{i}_t$.
The dependencies between random variables are illustrated in the probabilistic graphical model shown in Fig. \ref{fig:PGM}. 

After obtaining the sequence $U^*_t$, we apply only the first element, $u^*_t$, to the system and record a new state measurement $x_{t+1}$. Then, the dataset $\mathcal{D}_t$ in (\ref{eq:data_set_D}) is updated with the new data point $(\xi_t, x_{t+1})$, leading to an improved estimation of the dynamics, and the whole procedure will be repeated for the next time step, $t+1$. 

In the following sections, we introduce an inference method for real-time dynamics learning, formulate an AIF-based dual MPC scheme, and provide an approximate solution to the resulting SOC problem.

\section{Probabilistic Inference for Model Learning} \label{sec:Inf4Model}

This section presents a probabilistic inference approach for real-time identification of the unknown transition dynamics (\ref{eq:general_dynamics_transition}) and quantifying uncertainty based on the data collected in the dataset $\mathcal{D}_t$ defined in (\ref{eq:data_set_D}). We adopt a sparse semi-parametric GP model to provide probabilistic modelling for integrating into the model-based decision-making method introduced in the previous section. 
The semi-parametric GP model combines parametric and non-parametric components to capture both structured and unstructured uncertainties in the dynamics. By providing covariance estimates alongside the predicted mean, GPs can identify regions in the state-action space where prediction quality is poor due to limited data coverage or inherent stochasticity \cite{Rasmussen2006}. 
Moreover, to further reduce the computational cost of the GP, in addition to using a parametric part, we propose updating the GP online during the receding horizon control and applying a sparsification method to keep the dataset small.

\subsection{Semi-Parametric Gaussian Process}
At each real-time step $t$, we consider the following semi-parametric model, using the datasets $\mathcal{D}^{\xi}_t$ and $\mathcal{D}^{x_i}_t$ defined in (\ref{eq:data_set_D_split}):
\begin{align} 
    \hat{f}_{i}(\xi) &= \theta_{i}^{\top} \phi_{i}(\xi) + r_{i}, \label{eq:GP_f_decom}\\
    r_{i} &\sim \mathcal{GP}(0, \kappa_i(\xi,\xi')), \label{eq:GP_residual}
\end{align}
where $\theta_{i}^{\top} \phi_{i}(\xi)$ represents the linearly parameterized part of the model in which $\phi_{i}(\xi) \in \mathbb{R}^{n_{\theta_i}}$ is a feature transformation (e.g., a basis function) and assumed to be bounded, i.e., $\| \phi_{i}(\xi) \|_2 \le \bar{\phi}_{i}$; $r_{i}(\xi)$ denotes a zero-mean Gaussian process that models the residual. 
The parametric uncertainty $\theta_i \in \mathbb{R}^{n_{\theta_i}}$ is assumed to follow a Gaussian prior $p(\theta_i) = \mathcal{N}(\theta_i;\mu_{\theta_i,0}, \Sigma_{\theta_i,0})$.
$\kappa_i(\cdot,\cdot)$ denotes the covariance function, also known as a kernel. This work employs a squared exponential kernel defined as 
\begin{equation} \label{eq:GP_Kernel}
    \kappa_i(\xi,\xi') \coloneqq A_i \exp\left(-\tfrac{1}{2} \lVert \xi-\xi' \rVert^2_{W_i^{-1}} \right) + \sigma^2_{f_i} \delta_{ab}. 
\end{equation}
where $\delta_{ab}$ denotes the Kronecker delta, equal to one when $a = b$ and zero otherwise. The set $\psi_i \coloneqq \{ A_i, W_i, \sigma^2_{f_i} \}$ gathers tunable hyperparameters used to configure the kernel function \cite{Rasmussen2006, Pan2014, Pan2018, Deisenroth2015, Boedecker2014}. By integrating the prior distribution of $\theta_i$ in (\ref{eq:GP_f_decom}) and (\ref{eq:GP_residual}), we obtain a non-zero GP as follows: 
\begin{equation} \label{eq:GP_semiparametric}
    f_{i}(\xi) \sim \mathcal{GP}(\phi^{\top}_{i}(\xi)\mu_{\theta_i,0} , \tilde{\kappa}_i (\xi,\xi'))
\end{equation}
with a new kernel function 
\begin{equation}
    \tilde{\kappa}_i (\xi,\xi') \coloneqq \kappa_i(\xi,\xi')+\phi^{\top}_{i}(\xi)\Sigma_{\theta_i,0}\phi_{i}(\xi').
\end{equation}
In this article, we consider the following assumption.
Basically, the assumption guarantees that if the number of data points increases the error converges to zero. 
\begin{assumption} \label{assum:RKHS}
    Let for each $i \in \{1,\dots,n_x\}$, $f_i \in \mathcal{H}_{i}$, where $\mathcal{H}_{i}$ denotes the reproducing kernel Hilbert space (RKHS) associated with the kernel $\tilde{\kappa}_i$ and $\|\cdot\|_{\mathcal{H}_{i}}$ represents the induced norm of $\mathcal{H}_{i}$ \cite{Rasmussen2006, Hashimoto2025}.
\end{assumption}

The predictive distribution of the Gaussian Process for the $i^{\text{th}}$ element of the system at time $t$, given a single state-action sample $\xi^*_k$ or a batch of samples $\mathcal{D}^{\xi,*}_t = \{ \xi^{*,j}_k \mid j = 0, 1, \dots \}$, is given by:
\begin{equation}
    \label{eq:GP_pred_f}
    p^{\mathcal{GP}}(\hat{f}_{i,k+1}|\xi^*_{k}, \mathcal{D}^i_t) = \mathcal{N}(\hat{f}_{i,k+1}; \mu^{\mathcal{GP}}_{i,t} (\xi^*_k),\Sigma^{\mathcal{GP}}_{ii,t} (\xi^*_k)),
\end{equation}
where the mean and covariance functions, as adopted from \cite{Rasmussen2006}, are given by:
\begin{equation} \label{eq:pred_SemiGP}
    \begin{aligned}
        \mu^{\mathcal{GP}}_{i,t} (\xi^*) &=  \Phi_{i}^{*,\top}(\xi^*)\mu_{\theta_i,t}  + \mathcal{K}_{i,t}^{*,\top}(\xi^*) \mathcal{K}_{i,t}^{-1} (Y_{i,t} - \Phi_{i,t}^{\top}\mu_{\theta_i,t}), \\
        \Sigma^{\mathcal{GP}}_{ii,t} (\xi^*) &= \begin{multlined}[t]
        \mathcal{R}_{i,t}^{\top}(\xi^*)\Sigma_{\theta_i,t}\mathcal{R}_{i,t}(\xi^*)\\
        + \mathcal{K}_{i}^{**}(\xi^*)
        - \mathcal{K}_{i,t}^{*,\top}(\xi^*) \mathcal{K}_{i,t}^{-1} \mathcal{K}_{i,t}^{*}(\xi^*),
        \end{multlined}
\end{aligned}
\end{equation}
in which
\begin{subequations} \label{eq:GP_matrices}
    \begin{align}
        Y_{i,t} &\coloneqq \text{vec}(\mathcal{D}^{x_i}_t) \in \mathbb{R}^{|\mathcal{D}^i_{t}|}, \\
        \mathcal{K}_{i,t} &\coloneqq \kappa_{i}(\mathcal{D}_t^{\xi},\mathcal{D}_t^{\xi})\in \mathbb{R}^{|\mathcal{D}^i_{t}| \times |\mathcal{D}^i_{t}|},\\
        \mathcal{K}_{i,t}^{*}(\xi^*) &\coloneqq \kappa_{i}(\mathcal{D}_t^{\xi}, \mathcal{D}_t^{\xi,*})\in \mathbb{R}^{|\mathcal{D}^i_{t}| \times |\mathcal{D}^{\xi,*}_{t}|}, \\
        \mathcal{K}_{i}^{**}(\xi^*) &\coloneqq \kappa_{i}(\mathcal{D}_t^{\xi,*}, \mathcal{D}_t^{\xi,*})\in \mathbb{R}^{|\mathcal{D}^{\xi,*}_{t}| \times |\mathcal{D}^{\xi,*}_{t}|},\\
        \Phi_{i,t} & \coloneqq \phi_i(\mathcal{D}_t^{\xi})\in \mathbb{R}^{n_{\theta_i} \times |\mathcal{D}^i_{t}|},\\
        \Phi^{*}_{i} (\xi^*) & \coloneqq \phi_i(\mathcal{D}_t^{\xi,*}) \in \mathbb{R}^{n_{\theta_i} \times |\mathcal{D}^{\xi,*}_{t}|},\\
        \mathcal{R}_{i,t}(\xi^*) & \coloneqq  \Phi_{i,t} \mathcal{K}_{i,t}^{-1} \mathcal{K}_{i,t}^{*}(\xi^*) - \Phi^{*}_{i} (\xi^*),
    \end{align}
\end{subequations}
and the moments of $\theta_i$ are derived as follows:
\begin{subequations}\label{eq:postrior_params_GP}
    \begin{align}
        \mu_{\theta_i,t} &= \Sigma_{\theta_i,t}\Bigl(\Phi_{i,t} \mathcal{K}_{i,t}^{-1}Y_{i,t} + \Sigma_{\theta_i,0}^{-1}\mu_{\theta_i,0}\Bigr),  \label{eq:postrior_params_GP_mu}\\
        \Sigma_{\theta_i,t} &= \Bigl(\Phi_{i,t}\mathcal{K}_{i,t}^{-1} \Phi_{i,t}^{\top} + \Sigma_{\theta_i,0}^{-1}\Bigr)^{-1}.\label{eq:postrior_params_GP_sigma}
    \end{align}
\end{subequations}
\begin{remark}\label{rem:liklihood_theta}
    The expressions in (\ref{eq:postrior_params_GP}) correspond to the moments of a posterior distribution obtained when the likelihood is modeled as $p(\mathcal{D}^{i}_{t} \mid \theta_i) = \mathcal{N}(Y_{i,t} \mid \Phi_{i,t}^{\top} \theta_{i}, \mathcal{K}_{i,t})$.
\end{remark}

\begin{remark}
    The first terms in (\ref{eq:pred_SemiGP}) correspond to the parametric uncertainty, while the remaining terms represent the non-parametric component, which predicts the residual error between the observed data and the parametric model.
\end{remark}

\begin{lemma} \label{lem:GP_properties}
    If Assumption \ref{assum:RKHS} holds, then for a single state-action
    sample $\xi^*$, we have
    \begin{subequations}
        \begin{align}
            \text{I)} \quad &  \Sigma^{\mathcal{GP}}_{ii,t} (\xi^*) \le \bar{\Sigma}^{\mathcal{GP}}_{ii,t} \coloneqq A_i + \sigma^2_{f_i}+ \lambda_{\text{max}}[\Sigma_{\theta_i,0}] \bar{\phi}^2_{i}. \label{eq:GP_bound}\\
            \text{II)} \quad & \Sigma^{\mathcal{GP}}_{ii,t} (\xi^*) \rightarrow \sigma^2_{f_i} \text{ as } |\mathcal{D}^i_t| \rightarrow \infty. \label{eq:GP_consistency}
        \end{align}
    \end{subequations}
    where $\lambda_{\text{max}}[\cdot]$ returns the maximum eigenvalue of a matrix. 
\end{lemma}
\begin{proof}
Refer to Appendix \ref{app:lem1}. \qedhere
\end{proof}

\subsection{Real-Time Update of Semi-Parametric GP} 
Due to the receding horizon strategy, a new data point is observed at each time step. While incorporating this data improves prediction accuracy, recomputing the relevant matrices and functions in (\ref{eq:GP_matrices}) can become computationally expensive, especially when the cardinality of the dataset grows. 
To address this issue, we introduce a real-time update of the mean and covariance structure of the semi-parametric GP considered in this work. This approach extends the GP reparameterization techniques proposed in \cite{Csató2002, Chowdhary2015} to the semi-parametric GP setting introduced in the previous section. By defining the following parameters for the $i^{\text{th}}$ subsystem 
\begin{equation}\label{eq:def_recursive_param}
\begin{gathered}
\alpha_{i,t}\coloneqq \mathcal{K}_{i,t}^{-1}Y_{i,t}, \quad \beta_{i,t}\coloneqq\mathcal{K}_{i,t}^{-1}\Phi_{i,t}^{\top}, \quad \mathcal{C}_{i,t}\coloneqq\mathcal{K}_{i,t}^{-1}, 
\end{gathered}
\end{equation}
the moments of the predictive distribution of the semi-parametric GP in (\ref{eq:pred_SemiGP}) can be reformulated in a parameterized form as follows:
\begin{equation} \label{eq:GP_reparameterization}
    \begin{aligned}
        \mathcal{R}_{i,t}(\xi^*) &= \beta^{\top}_{i,t} \mathcal{K}_{i,t}^{*}(\xi^*) - \Phi^{*}_{i}(\xi^*), \\
        \mu^{\mathcal{GP}}_{i,t} (\xi^*) &=  \Phi_{i}^{*,\top}(\xi^*)\mu_{\theta_i,t}  + (\alpha_{i,t} - \beta_{i,t}\mu_{\theta_i,t})^{\top}\mathcal{K}_{i,t}^{*}(\xi^*),  \\
        \Sigma^{\mathcal{GP}}_{ii,t} (\xi^*) &= \begin{multlined}[t] \mathcal{R}_{i,t}^{\top}(\xi^*)\Sigma_{\theta_i,t}\mathcal{R}_{i,t}(\xi^*)\\
        +\mathcal{K}_{i}^{**}(\xi^*)-\mathcal{K}_{i,t}^{*,\top}(\xi^*) \mathcal{C}_{i,t} \mathcal{K}_{i,t}^{*}(\xi^*).
        \end{multlined}
    \end{aligned}
\end{equation}
When the next data pair $(\xi_t, x_{t+1})$ becomes available, the dataset $\mathcal{D}_t$ will be updated to $\mathcal{D}^i_{t+1} = \mathcal{D}^i_{t} \cup \{ (\xi_t, x_{i,t+1}) \}$. The parameters defined in (\ref{eq:def_recursive_param}) can then be recursively updated, as presented in the following lemma.

\begin{lemma}[Online Semi-Parametric GP] \label{lem:Online_SemiParametric_GP}
    The recursive update rules for $\alpha_{i,t}$, $\beta_{i,t}$, $\mathcal{C}_{i,t}$, $\mu_{\theta_i,t}$, and $\Sigma_{\theta_i,t}$ are given as follows: 
    \begin{subequations} \label{eq:GP_recursive_update}
        \begin{align}
            \alpha_{i,t+1} &= \mathcal{T} [\alpha_{i,t}] + a_{i,t}(\xi_t) S_{i,t}(\xi_t), \label{eq:GP_recursive_update_alpha}\\
            \beta_{i,t+1} &= \mathcal{T} [\beta_{i,t}] + b_{i,t}(\xi_t) S_{i,t}(\xi_t) \mathcal{R}_{i,t}^{\top}(\xi_t),
            \label{eq:GP_recursive_update_beta}\\
            \mathcal{C}_{i,t+1} &= \mathcal{U} [\mathcal{C}_{i,t}] + b_{i,t}(\xi_t) S_{i,t}(\xi_t) S_{i,t}^{\top}(\xi_t), \label{eq:GP_recursive_update_C}\\
            \mu_{\theta_i,t+1} &= \mu_{\theta_i,t} + \Lambda_{i,t}(\xi_t)\left(\mu^{\mathcal{GP}}_{i,t} (\xi_t) - x_{i,t+1} \right), \label{eq:GP_recursive_update_mu}\\
            \Sigma_{\theta_i,t+1} &= \left(\mathbb{I}_{n_{\theta_i}}- \Lambda_{i,t}(\xi_t)\mathcal{R}_{i,t}^{\top}(\xi_t) \right) \Sigma_{\theta_i,t}, \label{eq:GP_recursive_update_sigma}
        \end{align}
    \end{subequations}
    in which 
    \begin{subequations}
        \begin{align}
            b_{i,t}(\xi_t) &\coloneqq 1/(\mathcal{K}^{**}_{i,t}(\xi_{t}) - \mathcal{K}^{*,\top}_{i,t} (\xi_{t})\mathcal{C}_{i,t} \mathcal{K}^{*}_{i,t}(\xi_{t})),\\
            a_{i,t}(\xi_t) &\coloneqq b_{i,t}(\xi_t)\left(\alpha^{\top}_{i,t}\mathcal{K}^{*}_{i,t} (\xi_{t})-x_{i,t+1}\right),\\ 
            S_{i,t}(\xi_t) &\coloneqq \mathcal{T} [\mathcal{C}_{i,t} \mathcal{K}^{*}_{i,t}(\xi_{t})] - e,\\
            \Lambda_{i,t}(\xi_t) &\coloneqq \Sigma_{\theta_i,t} \mathcal{R}_{i,t} (\xi_t)\Sigma^{\mathcal{GP},-1}_{ii,t} (\xi_t), \label{eq:GP_recursive_update_Lambda}
        \end{align}
    \end{subequations}
    where $\phi_{i}(\xi_{t})$ represents the evaluation of the basis function at the new feature vector $\xi_{t}$, and $\mathcal{K}_{i,t}^{*}(\xi_{t})$ denotes the evaluation of the kernel function based on the dataset $\mathcal{D}^{\xi}_t$ collected up to time $t$. 
    The operators $\mathcal{T}[\cdot]$ and $\mathcal{U}[\cdot]$ extend a vector or matrix: $\mathcal{T}[\cdot]$ appends a zero row at the bottom, while $\mathcal{U}[\cdot]$ adds both a zero row at the bottom and a zero column on the right. The vector $e$ is the last standard unit vector with the same dimension as $|\mathcal{D}^i_{t+1}|$, with all entries equal to zero except for the last one, which is 1. 
    The initial values of the recursive parameters are determined based on the definition in (\ref{eq:def_recursive_param}) and the initial dataset $\mathcal{D}_0$, which may contain previously generated data and be initialized with the new data point $(\xi_0, x_1)$.
\end{lemma}
\begin{proof}
    Refer to Appendix \ref{app:lem2}. \qedhere
\end{proof}
The recursions in (\ref{eq:GP_recursive_update_alpha})–(\ref{eq:GP_recursive_update_sigma}) allow us to update the matrices defined in (\ref{eq:def_recursive_param}) without recalculating them from scratch at each time step.

\begin{remark}\label{rem:GP_recursions}
    The recursions in (\ref{eq:GP_recursive_update_alpha})–(\ref{eq:GP_recursive_update_C}) correspond to the nonparametric part of the semi-parametric GP, where the sizes of these matrices continue to grow as new data is observed, due to the operators $\mathcal{T}[\cdot]$ and $\mathcal{U}[\cdot]$. In contrast, the recursions in (\ref{eq:GP_recursive_update_mu}) and (\ref{eq:GP_recursive_update_sigma}) correspond to the parametric part, where the matrix sizes remain fixed. This parametric part efficiently absorbs new information into the existing matrices without adding additional computational cost related to the increase in cardinality of the dataset. 
\end{remark}

\subsection{Sparse Semi-Parametric GP}  
As noted in Remark \ref{rem:GP_recursions}, after observing new data and applying the recursions in (\ref{eq:GP_recursive_update_alpha})–(\ref{eq:GP_recursive_update_sigma}), the sizes of the matrices $\alpha_{i,t+1}$, $\beta_{i,t+1}$, and $\mathcal{C}_{i,t+1}$ continue to grow.
As a result, evaluating the mean and covariance of the GP in (\ref{eq:GP_reparameterization}) for a large dataset becomes computationally expensive. Therefore, after observing the new data, if the computational budget is exceeded, an existing data point in the dataset $\mathcal{D}^i_{t+1}$ corresponding to each subsystem must be removed. We explore the strategy adopted from \cite{Csató2002, Chowdhary2015} and extend it to the online semi-parametric GP introduced in the previous section. 
This method considers the Reproducing Kernel Hilbert Space (RKHS) norm with respect to the kernel function defined in (\ref{eq:GP_Kernel}) for the difference between the mean function of the current GP and alternative GPs, each missing a different data point, and then eliminates the data point with the lowest RKHS norm. 
The elimination score and deletion update equations for $\alpha_{i,t+1}$, $\beta_{i,t+1}$, and $\mathcal{C}_{i,t+1}$ are presented in the following lemma. 

\begin{lemma}[Sparse Semi-Parametric GP] \label{lem:Sparse_SemiParametric_GP}
A data point from the dataset $\mathcal{D}^i_{t+1}$ is deleted by selecting the index corresponding to the lowest elimination score based on the RKHS norm, which is given by
\begin{equation} \label{eq:elimination_score}
    \varepsilon_{i,t+1}(j)  \coloneqq  \frac{\mathcal{V}_{j}[|\alpha_{i,t+1} - \beta_{i,t+1}\mu_{\theta_i,t+1}|]}{\mathcal{V}_{j,j}[\mathcal{C}_{i,t+1}]},
\end{equation}
where $j=1,\dots,|\mathcal{D}^i_{t+1}|$, the operator $\mathcal{V}_{j}[\cdot]$ extracts the $j^\text{th}$ row of a matrix or vector, while $\mathcal{V}_{j,j}[\cdot]$ returns the matrix element at the $j^\text{th}$ row and the $j^\text{th}$ column.
Let $m$ be the data point selected for removal based on the score mentioned above. The deletion equations are then given by
\begin{equation} \label{eq:deletion_dataset}
        \mathcal{D}^{i}_{t+1} =  \mathcal{D}^{i}_{t+1} \backslash \{(\xi_{m-1}, x_{i,m})\},
\end{equation}

\begin{subequations} \label{eq:deletion_update}
    \begin{align}
        \alpha_{i,t+1} &= \mathcal{T}^{-}_{m} [\alpha_{i,t+1}] - \frac{\mathcal{T}^{-}_{m} [\mathcal{V}_{m}[\mathcal{C}_{i,t+1}]]^{\top} \mathcal{V}_{m}[\alpha_{i,t+1}]}{\mathcal{V}_{m,m}[\mathcal{C}_{i,t+1}]}, \label{eq:deletion_update_alpha} \\
        \beta_{i,t+1} &= \mathcal{T}^{-}_{m} [\beta_{i,t+1}] - \frac{\mathcal{T}^{-}_{m} [\mathcal{V}_{m}[\mathcal{C}_{i,t+1}]]^{\top} \mathcal{V}_{m}[\beta_{i,t+1}]}{\mathcal{V}_{m,m}[\mathcal{C}_{i,t+1}]}, \label{eq:deletion_update_beta} \\
        \mathcal{C}_{i,t+1} &= \mathcal{U}^{-}_{m,m} [\mathcal{C}_{i,t+1}] - \frac{\mathcal{T}^{-}_{m} [\mathcal{V}_{m}[\mathcal{C}_{i,t+1}]]^{\top} \mathcal{T}^{-}_{m} [\mathcal{V}_{m}[\mathcal{C}_{i,t+1}]]}{\mathcal{V}_{m,m}[\mathcal{C}_{i,t+1}]}, \label{eq:deletion_update_C}
    \end{align}
\end{subequations}
where the operators $\mathcal{T}^{-}_{m}[\cdot]$ and $\mathcal{U}^{-}_{m,m}[\cdot]$ reduce the dimension of a vector by removing its $m^\text{th}$ row and a matrix by removing its $m^\text{th}$ row and column, respectively. 
\end{lemma}
\begin{proof}
    Refer to Appendix \ref{app:lem3}. \qedhere
\end{proof}

\begin{remark}
    The elimination score (\ref{eq:elimination_score}) is proportional to both the contribution of each data point to the residual error between the parametric and non-parametric parts, and its share in the posterior covariance matrix, which indicates the GP's uncertainty about the dynamics. 
\end{remark}

\section{Active Inference for Dual Control} \label{sec:ActInf}

This section presents a dual control approach based on AIF. To design such a dual control, we leverage the POMDP formulation of the dual control problem, in which parametric uncertainties are treated as latent variables \cite{Mesbah2018}. In this work, we extend this perspective by treating the entire uncertain dynamics function, $\hat{f}_k$, including both its parametric and nonparametric components from the probabilistic model learning approach described in Sec. \ref{sec:Inf4Model}, as a latent variable. The observations $x_k$ then serve as outputs from the unknown environment over a finite horizon, as demonstrated in Fig. \ref{fig:PGM}. Finally, we apply the EFE minimization principle \cite{Kouw2024, Millidge2021}, whose underlying intent aligns with dual control: encouraging the agent to select actions that both reduce uncertainty about latent states and guide the system toward preferred outcomes.

The goal is to infer a variational distribution over both the latent variables $F_t$ and actions $U_t$ that optimally fits a generative model biased towards a task (see later) \cite{Kouw2024, Millidge2021}. 

Given the dataset $\mathcal{D}_t$ at each time step, we can infer the following posterior distribution for the latent states $\hat{f}_k$ based on the predictive distribution provided by the GP in (\ref{eq:pred_SemiGP}). 
\begin{equation} \label{eq:AIF_latent_posterior}
    p(F_{i,t}|\bar{Z}_t, \mathcal{D}^i_t) = \\ 
    \mathcal{N}(\hat{f}_{i,k+1}; \mu^{\mathcal{GP}}_{i,t} (\bar{Z}_t),\Sigma^{\mathcal{GP}}_{ii,t} (\bar{Z}_t)),
\end{equation}
where $\mu^{\mathcal{GP}}_{i,t} (\bar{Z}_t) \in \mathbb{R}^{H}$ and $\Sigma^{\mathcal{GP}}_{ii,t} (\bar{Z}_t) \in \mathbb{R}^{H\times H}$ are batch evaluations of the moments in (\ref{eq:pred_SemiGP}). 
This posterior distribution incorporates both the parametric and nonparametric parts of the estimation of dynamics functions as defined in (\ref{eq:GP_f_decom}), given the dataset $\mathcal{D}_t$. 
Furthermore, we consider the following noisy observation model based on the dynamic model in (\ref{eq:general_dynamics_transition}). 
\begin{equation}\label{eq:AIF_obs}
        p(x_{i,t+1}|\hat{f}_{i,t+1}) = \mathcal{N}(x_{i,k+1};\hat{f}_{i,k+1},\sigma^2_{f_i}),
\end{equation}
equivalently, over the horizon, we have
\begin{equation} \label{eq:AIF_gen_observation}
\begin{aligned}
    p(\bar{X}_{i,t+1}|F_{i,t}) &= \prod\nolimits_{k = t}^{t+H-1} p(x_{i,k+1}|\hat{f}_{i,k+1})  \\
    &= \mathcal{N}(\bar{X}_{i,t+1};F_{i,t},\sigma^2_{f_i}\mathbb{I}_{\text{H}}),
\end{aligned}
\end{equation}
where $\mathbb{I}_{\text{H}}$ is the identity matrix of size $H$.
The generative model for the AIF given the dataset $\mathcal{D}^i_t$ for the $i^{\text{th}}$ subsystem can be expressed using the following factorization \cite{Frigola2014}:
\begin{equation} \label{eq:AIF_gen_factorization}
\begin{multlined}
p(X_{i,t},F_{i,t}|U_t, \mathcal{D}^i_t) \\
\begin{aligned}
&\begin{multlined}
= p(x_{i,t})\prod\nolimits_{k = t}^{t+H-1} p(x_{i,k+1}|\hat{f}_{i,k+1}) \\ \times  p(\hat{f}_{i,k+1}|\hat{f}_{i,t+1:k}, \xi_{t:k},\mathcal{D}^i_t)
\end{multlined} \\
&= p(x_{i,t}) p(\bar{X}_{i,t+1}|F_{i,t}) p(\hat{f}_{i,t+1:t+H} | \xi_{t:t+H-1}, \mathcal{D}^i_t) \\
&= p(x_{i,t}) p(\bar{X}_{i,t+1}|F_{i,t}) p(F_{i,t}|\bar{Z}_t, \mathcal{D}^i_t).
\end{aligned}
\end{multlined}
\end{equation}
The marginal generative distribution is given by
\begin{equation} \label{eq:AIF_gen_x}
\begin{multlined}
p(\bar{X}_{i,t+1}|U_t , \mathcal{D}^i_t) \\
\begin{aligned}
&= \int p(\bar{X}_{i,t+1}|F_{i,t}) p(F_{i,t}|\bar{Z}_t, \mathcal{D}^i_t) \text{d}F_{i,t}\\
&= \mathcal{N}(\bar{X}_{i,t+1}; \mu^{\mathcal{GP}}_{i,t} (\bar{Z}_t),\Sigma^{\mathcal{GP}}_{ii,t} (\bar{Z}_t) + \sigma^2_{f_i}\mathbb{I}_{\text{H}}).
\end{aligned}
\end{multlined}
\end{equation}

The active inference problem seeks to minimize the approximation error between the variational distribution $q(F_t, U_t)$ and the biased generative model $p^*(Z_t, F_t)$. Within the setting of AIF, this approximation error is quantified by the following variational objective function known as the EFE functional \cite{Kouw2024, Millidge2021}.
\begin{equation}\label{eq:AIF_EFE}
        \mathcal{F}^{\text{EFE}}_{t} [q] \coloneqq \mathbb{E}_{q(Z_t,F_t)} 
        \left[ \ln \frac{q(U_t, F_t)}{p^*(Z_t, F_t)} \right]
\end{equation}
Using (\ref{eq:AIF_latent_posterior}) and (\ref{eq:AIF_gen_factorization}), we have the following decompositions: 
\begin{align} \label{eq:AIF_variational_decompositions}
    q(Z_t, F_t) &= p(X_t,F_t|U_t, \mathcal{D}_t) q(U_t), \\
    q(U_t, F_t) &= p(F_t|U_t, \mathcal{D}_t) q(U_t), \\
    p^*(Z_t, F_t) &= p(F_t|\bar{Z}_t, \mathcal{D}_t)p^*(Z_t).
\end{align}
Here, $q(U_t)$ denotes the approximate distribution over actions that we aim to infer and $p^*(Z_t)$ represents the \textit{goal prior distribution}, which encodes the notion of optimality based on the cumulative cost function (\ref{eq:cumulative_cost}) defined as: 
\begin{equation}
    p^*(Z_t) \propto \exp(-C_t(Z_t)). 
\end{equation} 

By substituting (\ref{eq:AIF_variational_decompositions}) into the EFE functional (\ref{eq:AIF_EFE}), one can obtain
\begin{align}
    \mathcal{F}^{\text{EFE}}_{t} [q] &= \mathbb{E}_{p(X_t,F_t|U_t,\mathcal{D}_t) q(U_t) } \Biggl[
    \ln \frac{p(F_t \mid U_t, \mathcal{D}_t) q(U_t)}{p(F_t|\bar{Z}_t , \mathcal{D}_t)p^*(Z_t)}
    \Biggr]
    \nonumber \\
    &=\mathbb{E}_{q(U_t)} \Bigl[ \ln \frac{q(U_t)}{\exp \bigl(-\mathcal{J}^{\text{EFE}}_{t} (U_t)\bigr)}
    \Bigr]. \nonumber
\end{align}
where $\mathcal{J}^{\text{EFE}}_{t} (U_t)$ is called \textit{EFE function} defined as:
\begin{equation}
    \mathcal{J}^{\text{EFE}}_{t} (U_t) \coloneqq
      \mathbb{E}_{p(X_t, F_t \mid U_t, \mathcal{D}_t)} \Biggl[
        \ln \frac{p(F_t \mid U_t , \mathcal{D}_t)}
                  {p(F_t|\bar{Z}_t , \mathcal{D}_t)p^*(Z_t)}
      \Biggr].
\end{equation}
Therefore, the optimal variational distribution for the inference problem can be represented by the following Boltzmann distribution:
\begin{equation}
        q^*(U_t) = \argmin_{q} \mathcal{F}^{\text{EFE}}_{t} [q] = \exp \bigl(-\mathcal{J}^{\text{EFE}}_{t} (U_t)\bigr). 
\end{equation}
It enables us to minimize the EFE function directly to find the optimal actions under AIF, i.e.,
\begin{equation}
        U_t^* = \argmax_{U_t} q^*(U_t) = \argmin_{U_t} \mathcal{J}^{\text{EFE}}_{t} (U_t). 
\end{equation}
To simplify the EFE function, we apply the Bayes' rule and the factorization outlined in (\ref{eq:AIF_gen_factorization}), one can have
\begin{equation} \label{eq:J_EFE}
\begin{multlined}
\mathcal{J}^{\text{EFE}}_{t} (U_t) \\
\resizebox{\linewidth}{!}{$
\begin{aligned}
&\begin{multlined}
= \mathbb{E}_{p(X_t,F_t \mid U_t, \mathcal{D}_t)} \left[ -\ln p^*(Z_t) \right] \\ + \mathbb{E}_{p(\bar{X}_{t+1},F_t \mid U_t, \mathcal{D}_t)} \left[ \ln \frac{p(F_t \mid U_t; \mathcal{D}_t)p(\bar{X}_{t+1} \mid U_t, \mathcal{D}_t)}{p(\bar{X}_{t+1},F_t \mid U_t, \mathcal{D}_t)} \right]
\end{multlined} \\
&= \mathbb{E}_{p(X_t \mid U_t , \mathcal{D}_t)} \left[ C_t(Z_t) \right] - \mathcal{I}(\bar{X}_{t+1};F_t \mid U_t , \mathcal{D}_t), 
\end{aligned}
$}
\end{multlined}
\end{equation}
where $\mathcal{I}(\bar{X}_{t+1};F_t \mid U_t , \mathcal{D}_t)$ denotes the mutual information (MI) of $\bar{X}_{t+1}$ and $F_t$ for fixed $U_t$ and $\mathcal{D}_t$ defined as:
\begin{multline}\label{eq:MI}
    \mathcal{I}(\bar{X}_{t+1};F_t\mid U_t, \mathcal{D}_t)
    \coloneqq \mathbb{E}_{p(\bar{X}_{t+1},F_t \mid U_t , \mathcal{D}_t)} \bigg[ \\ 
    \ln \frac{p(\bar{X}_{t+1}, F_t \mid U_t , \mathcal{D}_t)}{p(F_t \mid U_t , \mathcal{D}_t)p(\bar{X}_{t+1} \mid U_t , \mathcal{D}_t)} \bigg]. 
\end{multline}
The above MI can be equivalently expressed as:
\begin{equation*}\label{eq:MI_decomposition}
    \mathcal{I}(\bar{X}_{t+1};F_t \mid U_t, \mathcal{D}_t) = \entropy{\bar{X}_{t+1} \mid U_t , \mathcal{D}_t} - \entropy{\bar{X}_{t+1} \mid F_t}.
\end{equation*}
The required entropies can be computed using (\ref{eq:AIF_gen_observation}) and (\ref{eq:AIF_gen_x}) as given by
\begin{equation} \label{eq:entropy_obs}
    \begin{aligned}
     \entropy{\bar{X}_{t+1} \mid F_t} &= \sum\nolimits_{i=1}^{n_x} \entropy{\bar{X}_{i,t+1} \mid F_{i,t}} \\
     &\propto \tfrac{1}{2} \sum\nolimits_{i=1}^{n_x} \ln |\sigma^2_{f_i}\mathbb{I}_{\text{H}}| \\
     &= \tfrac{H}{2} \sum\nolimits_{i=1}^{n_x} \ln \sigma^2_{f_i},
\end{aligned}
\end{equation}
\begin{equation}\label{eq:entropy_pred}
    \begin{multlined} 
        \entropy{\bar{X}_{t+1} \mid U_t , \mathcal{D}_t} \\
        \begin{aligned}
            &= \mathbb{E}_{p(\bar{X}_t \mid U_t, \mathcal{D}_t)} \left[  \sum\nolimits_{i=1}^{n_x} \entropy{\bar{X}_{i,t+1} \mid U_t , \mathcal{D}^i_t}  \right] \\
            &\propto \mathbb{E}_{p(\bar{X}_t \mid U_t, \mathcal{D}_t)} \left[  \sum\nolimits_{i=1}^{n_x}  \tfrac{1}{2} \ln | \Sigma^{\mathcal{GP}}_{ii,t} (\bar{Z}_t) + \sigma^2_{f_i}\mathbb{I}_{\text{H}}| \right].
        \end{aligned}
    \end{multlined}
\end{equation}
The proportional relationships are held equal to the same additive constants.   
Therefore, the MI defined in (\ref{eq:MI}) is obtained as 
\begin{multline} \label{eq:MI2}
    \mathcal{I}(\bar{X}_{t+1};F_t \mid U_t, \mathcal{D}_t)\\ \propto \mathbb{E}_{p(\bar{X}_t \mid U_t, \mathcal{D}_t)} \left[  \sum\nolimits_{i=1}^{n_x}  \tfrac{1}{2} \ln \frac{| \Sigma^{\mathcal{GP}}_{ii,t} (\bar{Z}_t) + \sigma^2_{f_i}\mathbb{I}_{\text{H}}|}{|\sigma^2_{f_i}\mathbb{I}_{\text{H}}|} \right].
\end{multline} 
Substituting it into (\ref{eq:J_EFE}), we have
\begin{multline}
    \mathcal{J}^{\text{EFE}}_{t} (U_t) \propto \mathbb{E}_{p(X_t \mid U_t, \mathcal{D}_t)} \bigg[ \sum\nolimits_{k=t}^{t+H} c_k(\xi_k) \\
    - \tfrac{1}{2} \sum\nolimits_{i=1}^{n_x}  \ln \frac{\left| \Sigma^{\mathcal{GP}}_{ii,t} (\bar{Z}_t) + \sigma^2_{f_i}\mathbb{I}_{\text{H}}\right|}{\left|\sigma^2_{f_i}\mathbb{I}_{\text{H}}\right|} \bigg].
\end{multline}
The aforementioned objective function does not exhibit an optimal substructure and is therefore not amenable to dynamic programming or Bellman’s optimality principle, preventing an efficient solution.  
To circumvent this issue, we assume sequential evaluations of (\ref{eq:AIF_latent_posterior}) and (\ref{eq:AIF_gen_x}) instead of the batch evaluations using the following factorization: 
\begin{equation} \label{eq:AIF_GP_Dynamics}
    \begin{multlined} 
        p(\bar{X}_{i,t+1}|U_t , \mathcal{D}^i_t) \\
        \begin{aligned}
            &= \prod\nolimits_{k = t}^{t+H-1} p(x_{i,k+1}|\xi_k , \mathcal{D}^i_t) \\
            &= \prod\nolimits_{k = t}^{t+H-1} \mathcal{N}(x_{i,k+1}; \mu^{\mathcal{GP}}_{i,t} (\xi_k),\Sigma^{\mathcal{GP}}_{ii,t} (\xi_k) + \sigma^2_{f_i}).
        \end{aligned}
    \end{multlined}
\end{equation}
This posterior corresponds to the probabilistic model in which the dependencies between the latent variables $f_{k}$ are assumed to be ignored; i.e., the dashed arcs in Fig. \ref{fig:PGM} are eliminated. Then, the entropy (\ref{eq:entropy_pred}) simplifies to
\begin{multline*}
    \entropy{\bar{X}_{t+1} \mid U_t , \mathcal{D}_t} \\
    \propto \mathbb{E}_{p(\bar{X}_t \mid U_t, \mathcal{D}_t)} \left[\sum\nolimits_{i=1}^{n_x} \sum\nolimits_{k=t}^{t+H-1} \tfrac{1}{2} \ln ( \Sigma^{\mathcal{GP}}_{ii,t} (\xi_k) + \sigma^2_{f_i}) \right].
\end{multline*}
Using it to compute MI, one can obtain 
\begin{equation} \label{eq:AIF_objective}
    \mathcal{J}^{\text{EFE}}_{t} (U_t) \propto \mathbb{E}_{p(X_t \mid U_t, \mathcal{D}_t)} \left[ c_{t+H}(\xi_{t+H})  + \sum\nolimits_{k=t}^{t+H-1} c^{\text{EFE}}_k(\xi_k) \right]
\end{equation}
in which
\begin{align}
    c_k^{\text{EFE}}(\xi_k) &\coloneqq c_{k}(\xi_k) + c_k^{\text{EXP}}(\xi_k), \label{eq:cost_EFE}\\
    c_k^{\text{EXP}}(\xi_k) &\coloneqq - \tfrac{1}{2} \sum\nolimits_{i=1}^{n_x} \ln \left(1+ \tfrac{\Sigma^{\mathcal{GP}}_{ii,t} (\xi_k)}{\sigma^2_{f_i}}\right). \label{eq:cost_exploration}
\end{align}
The new objective function can be solved by applying Bellman's optimality principle. It is a modified version of the original objective (\ref{eq:SOC_objective}), augmented with an additive exploration term that encourages the controller to gather more informative data by maximizing information gain. 
Therefore, within AIF, at each time step, we use the probabilistic model learned via online GP introduced in Section \ref{sec:Inf4Model} to optimize the objective function (\ref{eq:AIF_objective}). 
This objective actively guides the controller toward regions where the GP model exhibits high covariance and thus high uncertainty, which simultaneously enhances both control performance and the learning of the dynamics. 
In the next section, we employ a probabilistic trajectory optimization method to numerically solve the SOC problem arising from AIF.

\begin{remark}
{
    The superposition of information cost and task-related cost, as in the cost function (\ref{eq:cost_EFE}), has already been proposed in prior work \cite{Alpcan2015, Alpcan2011, Le2021, Ladislav2014}. However, in this work, we provide a rationale for it based on AIF and consider the exact SOC problem. 
}
\end{remark}

\begin{remark}
To ensure that the cost function (\ref{eq:cost_EFE}) remains positive and well-conditioned under the upper bound (\ref{eq:GP_bound}) in Lemma \ref{lem:GP_properties}, we add the following constant to (\ref{eq:cost_EFE}):
\begin{equation}
    \bar{c}_k \coloneqq  \tfrac{1}{2} \sum\nolimits_{i=1}^{n_x} \ln \left(1+ \tfrac{\bar{\Sigma}^{\mathcal{GP}}_{ii,t}}{\sigma^2_{f_i}}\right).
\end{equation}
Furthermore, due to the consistency property in (\ref{eq:GP_consistency}), as more informative data are collected, $\Sigma^{\mathcal{GP}}_{ii,t} (\xi_k) \rightarrow \sigma^2_{f_i}$, leading the exploration cost (\ref{eq:cost_exploration}) to converge to a constant and AIF to reduce to the standard stochastic optimal control problem described in Section \ref{sec:Problem_Formulation}.
\end{remark}

\begin{remark}
    If we consider the settings and simplifying assumptions introduced in \cite{Kouw2024}, the problem reduces to learning the dynamics only by using the parametric model described in Section \ref{sec:Inf4Model}. In this case, the objective function (\ref{eq:AIF_objective}) collapses to the form of the objective function proposed by \cite{Kouw2024}, and the exploration cost (\ref{eq:cost_exploration}) reduces to:
    \begin{equation}\label{eq:cost_exploration_parametric}
        c_k^{\text{EXP}}(\xi_k) \coloneqq - \tfrac{1}{2} \sum\nolimits_{i=1}^{n_x} \ln \left(1+\phi^{\top}_{i}(\xi_k)\Sigma_{\theta_i,t}\phi_{i}(\xi_k)\right).
    \end{equation}
    Additionally, we address a limitation of the previous work \cite{Kouw2024}, where the resulting problem collapses into a deterministic optimization due to certain approximations. In this work, we account for the original SOC formulation. 
\end{remark}

\begin{remark}
In the context of system identification and adaptive control, a persistent excitation (PE) condition must be satisfied to ensure the convergence of the parameter estimation error. This condition checks the positive definiteness of the accumulated regressor covariance matrix, $\phi_i^\top(\xi_k)\phi_i(\xi_k)$, over the time horizon. In practice, this means collecting data that sufficiently excites the regressor function $\phi_i(\xi_k)$ \cite{Boyd1986}. 
In our framework, if we consider only parametric uncertainty and minimizing the exploration cost (\ref{eq:cost_exploration_parametric}), then it equivalently improves the satisfaction of a weighted PE condition by actively gathering data that excites the regressor in directions where the covariance matrix of the parameter estimates, $\Sigma_{\theta_i,t}$, has large eigenvalues, i.e., in the directions where the model is more uncertain.
Therefore, our framework actively enhances the satisfaction of the PE condition. 
\end{remark}

\section{Probabilistic Trajectory Optimization} \label{sec:Traj_Opt}

To address the stochastic dynamic optimization problem induced by AIF, various approximate methods proposed in the literature can be employed \cite{Theodorou2010_SDDP, Pan2014, Pan2018, Sutton2017}. In this section, we present a novel approximate solution method based on differential dynamic programming (DDP). This approach provides a model-based solution to a general finite-horizon stochastic optimal control problem with the objective functions (\ref{eq:SOC_objective}) or (\ref{eq:AIF_objective}), subject to the following stochastic dynamics with state–action-dependent noise:
\begin{equation} \label{eq:stochastic_dynamics_f}
    p(x_{k+1}|\xi_k) = \mathcal{N}\left(x_{k+1};f(\xi_k), \Sigma_f(\xi_k)\right),
\end{equation}
which generalizes the dynamics assumed in factorization (\ref{eq:AIF_GP_Dynamics}). 

We first review the general DP formulation as the theoretical foundation for stochastic optimal control. We then discuss an approximate formulation for the stochastic dynamics (\ref{eq:stochastic_dynamics_f}) using Fourier-Hermite DP (FHDP) and uncertainty propagation (UP), following the approach of \cite{Filabadi2025}.

\subsection{Differential Dynamic Programming (DDP)}
In the DP framework, the solution to the stochastic optimal control problem can be expressed in terms of the optimal state–action value function $Q^*_k:\mathbb{R}^{n_\xi}\mapsto \mathbb{R}^+$.
For each time step $k$ within the finite horizon $H$ (i.e., $k=t,t+1,\dots,t+H$), the optimal action, $u^*_k$, is determined by minimizing the optimal state-action value function $Q^*_k(\xi_k)$ \cite{Bellman1957, Sutton2017}.
\begin{equation} \label{eq:optimal_determinstic}
	u_k^*(x_k) = \argmin_{u_k} Q^*_k(\xi_k),  
\end{equation}
where $Q^*_k$ is defined by the following Bellman equation 
\begin{equation} \label{eq:Q_function} 
	Q^*_k(\xi_k) \coloneqq c_k({\xi}_k) + \expect{\mathcal{N}\left(x_{k+1};f(\xi_k), \Sigma_f(\xi_k)\right)}\left[ V^*_{k+1}({x}_{k+1})\right],  
\end{equation}
where $V^*_k:\mathbb{R}^{n_x}\mapsto\mathbb{R}^+$ denotes the optimal value function or (also known as the optimal cost-to-go). 
The value function satisfies the following backward Bellman optimality recursion:
\begin{equation} \label{eq:Bellman_Optimality}
    V^*_k(x_k) = \min_{u_k} Q^*_k(\xi_k),
\end{equation}
with the terminal condition $V_{t+H}(x_{t+H}) \coloneqq c_{t+H}(x_{t+H})$. 

In practice, solving this recursion explicitly is computationally intractable for general dynamics due to the curse of dimensionality. 
To tackle this challenge, a common approach is to approximate the state-action value function. 
For deterministic nonlinear systems, a widely adopted approach is Differential DP (DDP), which iteratively applies local quadratic approximations to the value functions, system dynamics, and stage cost function using second-order Taylor series and solves a sequence of quadratic subproblems within the standard DP framework to update the control actions \cite{jacobson1970differential, Weiwei2004, Pan2014, Pan2018, Boedecker2014, Hassan2023, Filabadi2024}.

In \cite{Theodorou2010_SDDP}, a modified version of DDP was proposed to deal with stochastic dynamics called Stochastic DDP by using second-order expansions of both mean and covariance of the dynamics. However, in such local approaches that rely on linearizing and quadratizing the nonlinear problem, the difference between deterministic and stochastic optimal control solutions is neglected. For instance, when the covariance of the dynamics is constant, the solution from \cite{Theodorou2010_SDDP} reduces to that of standard DDP. 
The application of DDP to Gaussian belief dynamics was first introduced in \cite{Pan2014, Pan2018} called Probabilistic DDP (PDDP), specifically for GP dynamic models.
In this approach, the Gaussian belief is obtained through forward uncertainty propagation, approximated by exact moment matching.
However, incorporating belief dynamics significantly increases the dimensionality of the problem.
Moreover, during policy evaluation and in the computation of belief–control costs, it is typically assumed that the control is independent of the belief.

To address these gaps, we exploit a statistical approximation of $Q^*_k(\xi_k)$ over the Gaussian local region obtained via forward uncertainty propagation, which captures the stochastic nature of the dynamics.
This enables us to apply the FHDP method \cite{Hassan2023} to solve the problem.
The next section provides a detailed description of this approach, referred to as Stochastic FHDP.

\subsection{Stochastic FHDP} \label{subsec:Stochastic_FHDP}

This section introduces Stochastic FHDP, an extension of the FHDP method \cite{Hassan2023} and designed to address stochastic optimal control problems. Previous approaches \cite{Theodorou2010_SDDP, Pan2014, Pan2018} rely on Taylor series expansions around a nominal deterministic trajectory, which do not fully capture the effects of the stochasticity.
In contrast, our method statistically approximates the optimal state-action value function $Q_k^*(\xi_k)$ within a Gaussian region that reflects the stochastic nature of the system. This region is defined by the nominal trajectory density, obtained through forward closed-loop simulation using nominal actions inherited by the PDDP framework. We approximate the nominal trajectory density at each time step $k$, $p(\xi_k)$, as a Gaussian distribution $\mathcal{N}(\xi_k; \mu_{\xi,k}, \Sigma_{\xi\xi,k})$.

In the following lemma adopted by \cite{sarkka2023, Gelb1974}, we apply a statistical approximation of $Q_k^*(\xi_k)$ over this Gaussian density. Following \cite{Gelb1974}, we achieve a more accurate representation of the state-action value function compared to deterministic Taylor series expansions. Also, this approach effectively embeds system's stochasticity into the approximation. 

\begin{lemma}[Statistical Approximation] \label{lem:SQ}
Consider the following surrogate quadratic approximation for $Q_k^*(\xi_k)$:
   \begin{equation} \label{eq:Quadratization_Q}
   	 \begin{aligned}
            Q_k^*(\xi_k) &\approx \hat{Q}_k^*(\xi_k)=\frac{1}{2} \begin{bmatrix}
				1 \\ \delta \xi_k
			\end{bmatrix}^\top \begin{bmatrix}
				2\hat{Q}^*_{0,k} & \hat{Q}_{\xi,k}^{*,\top}  \\
				\square & \hat{Q}_{\xi\xi,k}^*  \\
			\end{bmatrix}\begin{bmatrix}
				1 \\ \delta \xi_k
		\end{bmatrix},
    \end{aligned} 
   \end{equation}
where $\delta \xi_k \coloneqq \xi_k-\mu_{\xi,k}$ with the following matrix partitioning
\begin{equation}
        \hat{Q}_{\xi,k}^{*} \coloneqq 
	\begin{bmatrix}
		\hat{Q}^{*}_{x,k}   \\
		\hat{Q}^{*}_{u,k}   \\
	\end{bmatrix}, \quad
	\hat{Q}^{*}_{\xi\xi,k} \coloneqq 
	\begin{bmatrix}
		\hat{Q}^{*}_{xx,k} & \hat{Q}_{ux,k}^{*,\top}  \\
		\square & \hat{Q}^{*}_{uu,k} \\
	\end{bmatrix}.
\end{equation}
It minimizes the following mean squared error
\begin{equation} \label{eq:MSE_SQ}
    \expect{\mathcal{N}(\xi_k; \mu_{\xi,k}, \Sigma_{\xi\xi,k})}\left[ \left(Q_k^*(\xi_k) - \hat{Q}_k^*(\xi_k) \right)^2 \right],
\end{equation}
when 
\begin{equation} \label{eq:SQ_coef}
    \begin{aligned}
        \hat{Q}^*_{0,k} &\coloneqq \expect{\mathcal{N}(\xi_k; \mu_{\xi,k}, \Sigma_{\xi\xi,k})} \left[ Q_{k}^*(\xi_k) \right] - \tfrac{1}{2} \trace \left( \hat{Q}_{\xi\xi,k}^{*} \Sigma_{\xi\xi,k}\right),\\
        \hat{Q}_{\xi,k}^{*} &\coloneqq \expect{\mathcal{N}(\xi_k; \mu_{\xi,k}, \Sigma_{\xi\xi,k})} \left[ Q_{\xi,k}^*(\xi_k) \right],\\
        \hat{Q}_{\xi\xi,k}^{*} &\coloneqq \expect{\mathcal{N}(\xi_k; \mu_{\xi,k}, \Sigma_{\xi\xi,k})} \left[ Q_{\xi\xi,k}^*(\xi_k) \right],
    \end{aligned}
\end{equation}
where $ Q_{\xi,k}^*(\xi_k)$ and $Q_{\xi\xi,k}^*(\xi_k)$ are first and second order derivatives of $Q_k^*(\xi_k)$ with respect to $\xi_k$. 
\end{lemma}

\begin{proof}
    The result follows from setting the derivatives of the mean squared error in (\ref{eq:MSE_SQ}) with respect to the coefficients $\hat{Q}^*_{0,k}$, $\hat{Q}_{\xi,k}^{*}$, and $\hat{Q}_{\xi\xi,k}^{*}$ to zero, and then applying Stein's lemma to the resulting expressions \cite{sarkka2023, Gelb1974}. \qedhere
\end{proof}

The Gaussian region ($\mu_{\xi,k}, \Sigma_{\xi\xi,k}$) is the quantified uncertainty from the forward pass. By leveraging the approximation introduced in Lemma \ref{lem:SQ}, regions with higher uncertainty yield a more generalized approximation of $Q^*_k(\xi_k)$. Consequently, our approach dynamically adjusts the scale of the approximation’s generalization in direct proportion to the uncertainty level.

\begin{remark} \label{rem:SQ_FHS}
    According to the definition of the Fourier-Hermite series \cite{Hassan2023}, the surrogate quadratic expansion (\ref{eq:Quadratization_Q}), with coefficients specified in (\ref{eq:SQ_coef}), is equivalent to the second-order Fourier-Hermite series of $Q_k^*(\xi_k)$ (see Lemma 2 in \cite{Hassan2023} for more details).
\end{remark}

Remark \ref{rem:SQ_FHS} allows the coefficients in (\ref{eq:SQ_coef}) to be computed using the derivative-free procedure of FHDP \cite{Hassan2023}, which employs sigma-point methods for numerical integration. By combining uncertainty propagation (UP) in the forward pass and the approximation of Lemma \ref{lem:SQ} with the coefficient evaluation method of \cite{Hassan2023} in the backward pass, we develop an uncertainty-aware planning framework for stochastic environments. In this framework, the covariance that defines the Gaussian region for approximation adapts automatically rather than remaining fixed as in \cite{Hassan2023}. 
The framework iterates forward and backward passes until convergence, as detailed below.

\subsubsection{Backward Pass}

Following \cite{Hassan2023, Filabadi2024, Filabadi2025}, by apply a sigma-point method, the coefficients of $\hat{Q}^*_k(\xi_k)$ in (\ref{eq:SQ_coef}) can approximately be computed as follows:
\begin{equation}\label{eq:FH_Q_1}
    \begin{aligned}
        \hat{Q}_{\xi,k}^* &= \sqrt{\Sigma}_{\xi\xi,k}^{-1}\tilde{Q}_{\xi,k}^*,\\
        \hat{Q}_{\xi\xi,k}^* &= \sqrt{\Sigma}_{\xi\xi,k}^{-1}\tilde{Q}_{\xi\xi,k}^*\sqrt{\Sigma}_{\xi\xi,k}^{-1},
    \end{aligned}
\end{equation}
in which 
\begin{equation}\label{eq:FH_Q_3}
    \begin{aligned}
        \tilde{Q}_{\xi,k}^* &= \sum\nolimits_{n=1}^{N_\xi} w^{\xi}_n Q_k^*\left(\mu_{\xi,k}+\sqrt{\Sigma}_{\xi\xi,k} \epsilon^{\xi}_n\right) \epsilon^{\xi}_n, \\
        \tilde{Q}_{\xi\xi,k}^* &= \sum\nolimits_{n=1}^{N_\xi} w^{\xi}_n Q_k^*\left(\mu_{\xi,k}+\sqrt{\Sigma}_{\xi\xi,k} \epsilon^{\xi}_n\right) \left(\epsilon^{\xi}_n\epsilon^{\xi\top}_n -\matrixstyle{I}\right),
    \end{aligned}
\end{equation}
where $\sqrt{\Sigma}_{\xi\xi,k}$ denotes the Cholesky factor or another square root of the covariance matrix $\Sigma_{\xi\xi,k}$. $N_{\xi}$ denotes the number of sigma points used for numerical integration. The set of $n_{\xi}$-dimensional unit sigma points, 
$\{\epsilon^{\xi}_n \in \mathbb{R}^{n_{\xi}}\}_{n=1}^{N_{\xi}}$,  is deterministically chosen to capture the key statistical properties of the distribution, with $\{w^{\xi}_n \in \mathbb{R}\}_{n=1}^{N_\xi}$ representing their associated scalar weights. Two primary methods for selecting these sigma points based on Gauss–Hermite and spherical cubature integration techniques are detailed in \cite{sarkka2023}.
We also consider the following quadratic form of the optimal value function $V_t^*(x_k)$:
   \begin{equation} \label{eq:Quadratization_V}
            V_t^*(x_k) \approx \hat{V}_k^*(x_k)=\frac{1}{2} \begin{bmatrix}
				1 \\ \delta x_k
			\end{bmatrix}^\top \begin{bmatrix}
				\hat{V}^*_{0,k} & \hat{V}_{x,k}^{*,\top}  \\
				\square & \hat{V}_{xx,k}^*  \\
			\end{bmatrix}\begin{bmatrix}
				1 \\ \delta x_k
		\end{bmatrix}.
   \end{equation}
To evaluate $Q^*_k(\xi_k)$ at sigma points for constructing (\ref{eq:FH_Q_3}), we employ the backward recursion in (\ref{eq:Q_function}) as follows:
\begin{equation}\label{eq:Bellman_Q}
    \resizebox{\linewidth}{!}{$
            \begin{aligned}
                Q^*_{k}(\xi_k)&\approx c_k({\xi}_k) + \expect{\mathcal{N}\left(x_{k+1};f(\xi_k), \Sigma_f(\xi_k)\right)}\left[ \hat{V}^*_{k+1}({x}_{k+1})\right]\\
                &\begin{multlined}
                    = c_k({\xi}_k) + \tfrac{1}{2}\trace\left(\Sigma_{f}(\xi_k)\hat{V}_{xx,k+1}^*\right) +\hat{V}_{x,k+1}^{*,\top} \Delta \mu_{k+1} (\xi_k) \\ +\tfrac{1}{2}\Delta \mu_{k+1}^{\top}(\xi_k) \hat{V}_{xx,k+1}^* \Delta \mu_{k+1} (\xi_k) + \hat{V}_{0,k+1}^{*},
                \end{multlined}
            \end{aligned}
    $}
\end{equation}
where $\Delta \mu_{k+1} (\xi_k) \coloneqq f(\xi_k) - \mu_{x,k+1}$. 
The aforementioned backward pass initiates from a second-order approximation of $V^*_{t+H}(x_{t+H}) = c_{t+H}(x_{t+H})$. Reapplying the approximation from Lemma \ref{lem:SQ} on $c_{t+H}(x_{t+H})$ together with the sigma-point evaluation methods of \cite{Hassan2023, Filabadi2024, Filabadi2025} then yields
\begin{equation}
    \begin{aligned} \label{eq:Terminal_Value_Coef}
        \hat{V}_{x,t+H}^* &= \sqrt{\Sigma}_{xx,t+H}^{-1} \tilde{V}_{x,t+H}^*, \\
        \hat{V}_{xx,t+H}^* &= \sqrt{\Sigma}_{xx,t+H}^{-1}\tilde{V}_{xx,t+H}^{*}\sqrt{\Sigma}_{xx,t+H}^{-1},
    \end{aligned}
\end{equation}
in which
\begin{equation}
    \begin{aligned} \label{eq:Terminal_Value_FH}
        \tilde{V}_{x,t+H}^* &= \sum_{n=1}^{N_x} w^x_n  c_{t+H}(\mu_{x,t+H}+\sqrt{\Sigma}_{xx,t+H}\epsilon^x_n) \epsilon^x_n,\\
        \tilde{V}_{xx,t+H}^* &= \sum_{n=1}^{N_x} w^x_n  c_{t+H}(\mu_{x,t+H}+\sqrt{\Sigma}_{xx,t+H} \epsilon^x_n) (\epsilon^x_n {\epsilon^{x\top}_n} -\matrixstyle{I}),
    \end{aligned}
\end{equation}
where $\{\epsilon^{x}_n \in \mathbb{R}^{n_{x}}\}_{n=1}^{N_x}$ and $\{w^{x}_n \in \mathbb{R} \}_{n=1}^{N_x}$ denote the set of $n_x$-dimensional sigma points and their corresponding scalar weights, respectively.

Once the coefficients of (\ref{eq:FH_Q_1}) are evaluated, the optimal deviation for the control action can be approximated in the following affine form that minimizes the quadratic approximation of $Q^*_k(\xi_k)$:
\begin{equation} \label{eq:Linearization_action}
    \delta u^*_k = k_k^*+\matrixstyle{K}_k^* \delta x_k,
\end{equation}
where $k^*_k = -\hat{Q}_{uu,k}^{*,-1} \hat{Q}_{u,k}^*$ and 
$\matrixstyle{K}^*_k = -\hat{Q}_{uu,k}^{*,-1} \hat{Q}_{ux,k}^*$. Consequently, following (\ref{eq:Bellman_Optimality}), the coefficients of $V^*(x_k)$ are updated by substituting (\ref{eq:Linearization_action}) into (\ref{eq:Quadratization_Q}) as follows: 
\begin{equation} \label{eq:coef_V}
    \begin{aligned}
        \hat{V}_{x,k}^* &= \hat{Q}_{x,k}^* - \matrixstyle{K}_k^{*\top} \hat{Q}_{uu,k}^* k^*_k, \\
        \hat{V}_{xx,k}^* &= \hat{Q}_{xx,k}^*  - \matrixstyle{K}_k^{*\top} \hat{Q}_{uu,k}^* \matrixstyle{K}^*_k. 
    \end{aligned}
\end{equation}
Note that the constant coefficients in (\ref{eq:Quadratization_Q}) and (\ref{eq:Quadratization_V}) are idle terms, so we omit their evaluations for brevity.

\begin{remark}
    Control constraints such as $u_{\text{min}} \le u^*_k \le u_{\text{max}}$ can be taken into account in various ways in our algorithm. In particular, the approaches proposed in \cite{Tassa2014, Pan2018} can be exactly incorporated. 
\end{remark}

\subsubsection{Forward Pass}
\begin{figure}[!h]
  \centering
  \includegraphics[width=1\columnwidth]{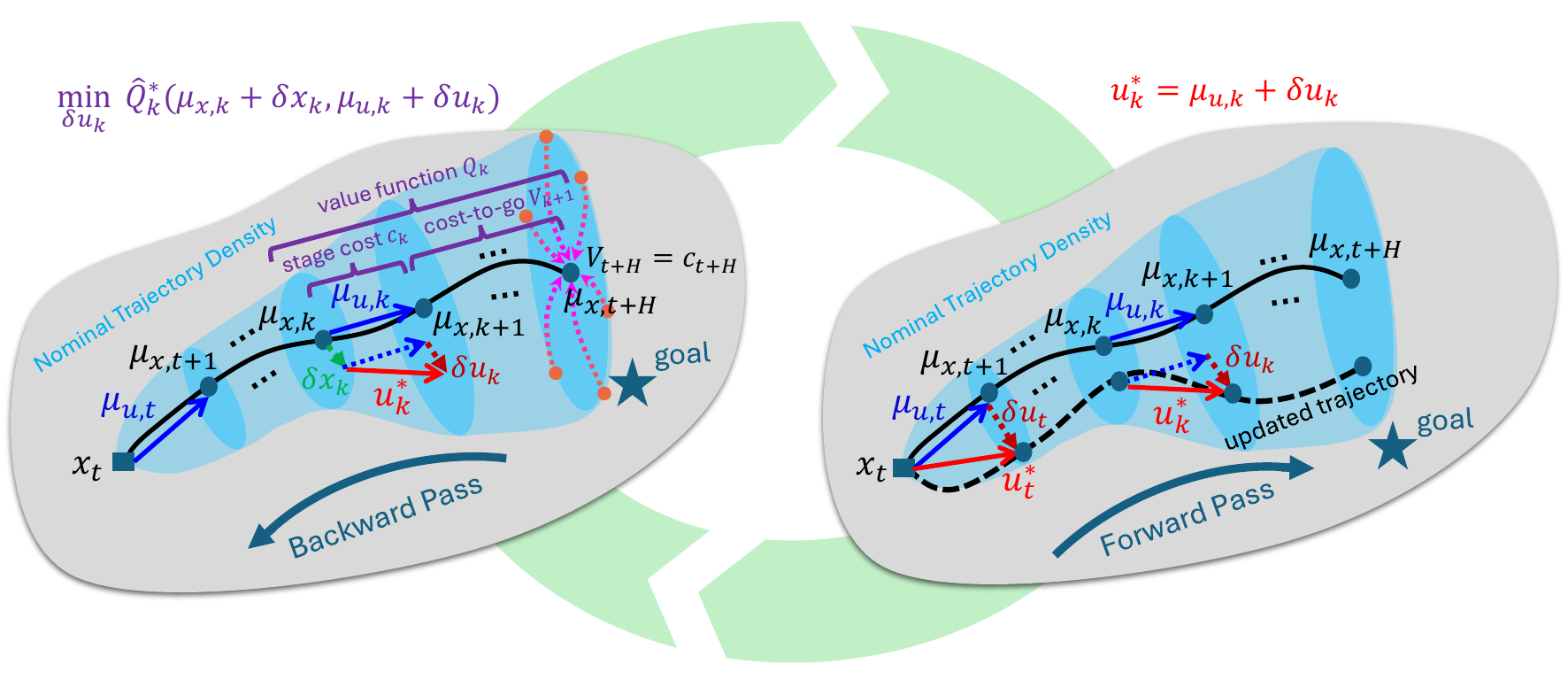} \caption{\small{
  Schematic representation of stochastic FHDP iterations: It incorporates trajectory densities (blue disks) through mean and variance representations $(\mu_{\xi,k}, \Sigma_{\xi\xi,k})$. The backward and forward passes are performed over probability distributions. 
  }}\label{fig:GP-FHDP}
\end{figure}

In the forward pass, the updated actions obtain from the backward pass are applied to propagate the system forward, resulting in an updated density $p(X_t^*, U_t^*)$ for the state-action trajectories. This density is assumed to follow a normal distribution $\mathcal{N}(\xi^*_k ; \mu^*_{\xi,k}, \Sigma^*_{\xi\xi,k})$ for each $k = t, \dots, t+H$. To that ends, the following forward uncertainty propagation recursion is computed, starting from the fixed initial distribution $\mathcal{N}(x_t; \mu_{x,t}, \Sigma_{xx,t})$. 
\begin{subequations} \label{eq:Forward_UP}
    \begin{align}
        p(\xi^*_k)   &= p(u^*_k|x^*_k) p(x^*_k) \label{eq:Forward_UP_1}, \\
        p(x^*_{k+1}) &= \int p(x^*_{k+1}|\xi^*_k) p(\xi^*_k) d\xi^*_k, \label{eq:Forward_UP_2}
    \end{align}
\end{subequations}
where $p(u^*_k | x^*_k)$ denotes the updated action policy, which is given by 
\begin{equation} \label{eq:updated_action}
        \delta\left(u^*_k = \mu_{u,k} + k^*_k + \matrixstyle{K}^*_k (x^*_k - \mu_{x,k})\right),
\end{equation}
where $\delta(\cdot)$ indicates the Dirac delta distribution, and $\mu_{x,k}$ and $\mu_{u,k}$ are determined by the nominal trajectory densities $p(\xi_k)$. Therefore, the statistical moments of (\ref{eq:Forward_UP_1}) can be computed as follows: 
\begin{equation} \label{eq:Forward_UP_Cont}
  \begin{aligned}
  	  \Sigma^*_{\xi\xi,k} &\coloneqq \begin{bmatrix}              				
      \Sigma^*_{xx,k} & \Sigma^*_{xu,k}  \\
             		\square & \Sigma^*_{uu,k} \\
                			    \end{bmatrix} = \begin{bmatrix}              				\Sigma^*_{xx,k} & \Sigma^*_{xx,k} \matrixstyle{K}_k^{*,\top}    \\
             		\square & \matrixstyle{K}^*_k \Sigma^*_{xx,k} \matrixstyle{K}_k^{*,\top} \\
                			    \end{bmatrix}, \\
          	 \mu^*_{\xi,k} &\coloneqq \begin{bmatrix}
  		\mu^*_{x,k}  \\
  		\mu^*_{u,k}  \\
  	\end{bmatrix} =  \begin{bmatrix}
  		\mu^*_{x,k}  \\
  		\mu_{u,k} + k^*_k + \matrixstyle{K}^*_k (\mu^*_{x,k}- \mu_{x,k})  \\
  	\end{bmatrix}. 
  \end{aligned}
\end{equation}
Finally, to compute the moments of the predictive distribution in (\ref{eq:Forward_UP_2}), we apply the exact Gaussian moment matching technique \cite{Pan2014, Pan2018, Deisenroth2015}, given by
\begin{equation} \label{eq:general_moment_matching}
    \begin{aligned}
              \mu^*_{x,k+1} &= \expect{\mathcal{N}(\xi^*_k;\mu^*_{\xi,k}, {\Sigma}^*_{\xi\xi,k})}  \left[ f(\xi^*_k) \right], \\
              \Sigma^*_{xx,k+1} &= \expect{\mathcal{N}(\xi^*_k;\mu^*_{\xi,k}, {\Sigma}^*_{\xi\xi,k})}  \left[ f(\xi^*_k)f^{\top}(\xi^*_k) \right]\\
              &+ \expect{\mathcal{N}(\xi^*_k;\mu^*_{\xi,k}, {\Sigma}^*_{\xi\xi,k})}  \left[ \Sigma_f(\xi^*_k) \right] - \mu^*_{x,k+1} \mu^{*,\top}_{x,k+1}.
    \end{aligned}
\end{equation}
The expectations can be computed using a $n_{\xi}$-dimensional sigma-point method. 
\begin{remark}\label{rem:sigma_points_forward_backward}
    Due to the fixed and deterministic generation of the unit sigma points, the evaluations of the dynamics at sigma points can be reused from the computations in (\ref{eq:FH_Q_3}) during the backward pass for enhancing computational efficiency or vice versa.
\end{remark}
In the next iteration of the algorithm the updated state-action trajectory density will be used as a nominal density.
Fig. \ref{fig:GP-FHDP} presents a schematic illustration of the stochastic FHDP algorithm, which contains repeated cycles of backward and forward passes until the optimal actions converge. 

\section{Dual Adaptive MPC Relying on AIF} \label{sec:algorithm}

Based on the previous sections, we propose a probabilistic MPC framework leveraging AIF to address the problem in Sec. \ref{sec:Problem_Formulation}. At each time step, backward and forward iterations from Sec. \ref{subsec:Stochastic_FHDP} are performed on the model learned via Sec. \ref{sec:Inf4Model} using the designed objective function from Sec. \ref{sec:ActInf}. Then, the first optimized control action is applied to the system, new data is collected, and the model is updated for the next time step. Some considerations and details of the algorithm’s pseudocode are explained below.

\subsubsection{Algorithm Overview}
In summary, the proposed algorithm, GPFHDP-AIF, is presented in Algorithm \ref{Alg:GPFHDP_AIF}, and its associated functions are described in Subroutine \ref{Alg:GPFHDP_AIF_Subroutines}. 
In the forward pass of probabilistic trajectory optimization (Sec. \ref{subsec:Stochastic_FHDP}), uncertainty is propagated through the dynamics to obtain the closed-loop trajectory density using (\ref{eq:general_moment_matching}). For a semi-parametric GP model at time $t$, the $i^{\text{th}}$ element of $f(\xi_k)$ and $\Sigma_f(\xi_k)$ is replaced by $\mu^{\mathcal{GP}}_{i,t} (\xi_k)$ and $\Sigma^{\mathcal{GP}}_{ii,t} (\xi_k)$, as in (\ref{eq:GP_reparameterization}). The cost function $c_k$ in (\ref{eq:Bellman_Q}) is replaced with $c^{\text{EFE}}_k$ (\ref{eq:cost_EFE}). Additionally, to better manage the trade-off between exploration and exploitation, we introduce a hyperparameter, $\gamma$, into the cost function (\ref{eq:cost_EFE}), resulting in $c_k^{\text{EFE}}(\xi_k) =  c_{k}(\xi_k) + \gamma c_k^{\text{EXP}}(\xi_k)$. 

As the algorithm iterates and collects more informative data, the uncertainty of the GP gradually decreases, ultimately allowing the GP to closely approximate the true dynamics of the system and allowing for more reliable decision-making.

\begin{algorithm}[!h] \small{
    \footnotesize 
	\caption{GPFHDP-AIF}
	\label{Alg:GPFHDP_AIF}
    \SetKwFunction{BackwardPass}{BackwardPass}
    \SetKwFunction{ForwardPass}{ForwardPass}
    \SetKwFunction{GPUpdate}{GP-Update}
	\KwIn{$H$, $c_{t:t+H}$, $\gamma$, $\mu_{x,0}$, $\Sigma_{xx,0}$, $\{\epsilon^{\xi}_n, w^{\xi}_n\}_{n=1}^{N_\xi}$, $\{\epsilon^{x}_n, w^{x}_n\}_{n=1}^{N_x}$, $\mu_{\theta_i,0}$, $\Sigma_{\theta_i,0}$, and ${\mathcal{D}_0}$.
	}


    Initialize $\alpha_{i,0}, \beta_{i,0}, \mathcal{C}_{i,0}$ using ${\mathcal{D}_0}$; \\
    
        $\mathcal{GP}_0 = \{ \alpha_{i,0}, \beta_{i,0}, \mathcal{C}_{i,0}, \mu_{\theta_i,0}, \Sigma_{\theta_i,0}\}_{i=1}^{n_x}$; \tcp*[f]{GP state} \\
    \For(\tcp*[f]{Real-Time Loop}){$t$}{ 
     Start from $\mu_{x,t}$ and $\Sigma_{xx,t}$; \\
     Initialize the nominal density $\{\mu_{\xi,t:t+H}, \Sigma_{\xi\xi, t:t+H}\}$; \\
	\Repeat{convergence}{
        $\{K^*_{t:t+H-1}, k^*_{t:t+H-1}\}$ = \BackwardPass{$\mu_{\xi,t:t+H}, \Sigma_{\xi\xi, t:t+H}, \mathcal{GP}_t$}; \\
        $\mu_{\xi,t:t+H}, \Sigma_{\xi \xi, t:t+H}$ = \ForwardPass{$\mu_{\xi,t:t+H}, K^*_{t:t+H-1}, k^*_{t:t+H-1}, \mathcal{GP}_t$}; \\ 
	}
    Compute $u^*_{t}$ using (\ref{eq:updated_action}); \\
    Apply $u^*_{t}$ on the dynamic system and observe $(\xi_{t},x_{t+1})$;  \\
    $\mathcal{GP}_{t+1}$ = \GPUpdate{$\xi_{t},x_{t+1}, \mathcal{GP}_{t}$}; \\
    $\{\mu_{x,t+1}, \Sigma_{xx,t+1}\} \leftarrow \{ \mu_{t+1}^{\mathcal{GP}}(\xi_{t}) , \Sigma_{t+1}^{\mathcal{GP}}(\xi_{t})\}$; \\
    }}
\end{algorithm}

\renewcommand{\algorithmcfname}{Subroutine}
\setcounter{algocf}{0} 
\begin{algorithm}[!h] \small{
    \footnotesize 
    \SetKwFunction{BackwardPass}{BackwardPass}
	\caption{GPFHDP-AIF Subroutines}
	\label{Alg:GPFHDP_AIF_Subroutines}
    \SetKwProg{myproc}{Function}{}{}
    \myproc{\BackwardPass{$\mu_{\xi,t:t+H}, \Sigma_{\xi\xi, t:t+H}, \mathcal{GP}_{t}$}}{
      Evaluate $\hat{V}_{x,t+H}^*$ and $\hat{V}_{xx,t+H}^*$ using (\ref{eq:Terminal_Value_Coef})-(\ref{eq:Terminal_Value_FH});\\
        \For(\tcp*[f]{Backward Recursion}){$k = t+H-1:t$}
		{
			Evaluate $Q^*_k$ at $\mu_{\xi,k}+\sqrt{\Sigma}_{\xi\xi,k} \epsilon^{\xi}_n$ using (\ref{eq:Bellman_Q}); \\
			Evaluate $\hat{Q}_{\xi,k}^*$ and $\hat{Q}_{\xi \xi,k}^*$ using (\ref{eq:FH_Q_1})-(\ref{eq:FH_Q_3}); \\
			Obtain $k^*_k$ and $K^*_k$ using (\ref{eq:Linearization_action}); \\ 
			Compute $\hat{V}_{x,k}^*$ and $\hat{V}_{xx,k}^*$ using (\ref{eq:coef_V}); \\
		}
      \KwRet $k^*_{t:t+H-1}, K^*_{t:t+H-1}$\;}

    \SetKwFunction{ForwardPass}{ForwardPass}
    \SetKwProg{myproc}{Function}{}{}
      \myproc{\ForwardPass{$\mu_{\xi,t:t+H}, K^*_{t:t+H-1}, k^*_{t:t+H-1}, \mathcal{GP}_{t}$}}{
        Start from $\mu_{x,t}$ and $\Sigma_{xx,t}$; \\
		\For(\tcp*[f]{Forward Recursion}){$k = t:t+H-1$}
		{
            Evaluate $\mu^*_{\xi,k+1}$ and $\Sigma^*_{\xi \xi,k+1}$ using (\ref{eq:Forward_UP_Cont})-(\ref{eq:general_moment_matching}); \\
		}
      \KwRet $\mu^*_{\xi,t:t+H}$ and $\Sigma^*_{\xi \xi,t:t+H}$\;}
      
      \SetKwFunction{GPUpdate}{GP-Update}
    \SetKwProg{myproc}{Function}{}{}
      \myproc{\GPUpdate{$\xi_t, x_{t+1}, \mathcal{GP}_{t}$}}{
        $\mathcal{D}_{t+1} \leftarrow \mathcal{D}_{t} \cup \{(\xi_{t},x_{t+1})\}$; \\
        \For(){$i = 1 : n_x$}
		{
            Update $\alpha_{i,t}, \beta_{i,t}, \mathcal{C}_{i,t}, \mu_{\theta_i,t}, \Sigma_{\theta_i,t}$ using (\ref{eq:GP_recursive_update}); \\
            Compute the elimination score $\varepsilon_{i,t+1}(j)$ using (\ref{eq:elimination_score});
            $m_i = \argmin_{j} \varepsilon_{i,t+1}(j)$; \\
            $\mathcal{D}^{i}_{t+1} =  \mathcal{D}^{i}_{t+1} \backslash \{(\xi_{m_i-1}, x_{i,m_i})\}$; \\
            Apply the update rule (\ref{eq:deletion_update}); \\
        }
        $\mathcal{GP}_{t+1} = \{ \alpha_{i,t+1}, \beta_{i,t+1}, \mathcal{C}_{i,t+1}, \mu_{\theta_i,t+1}, \Sigma_{\theta_i,t+1} \}_{i=1}^{n_x} $; \\
      \KwRet $\mathcal{GP}_{t+1}$}
      }
\end{algorithm}

\subsubsection{Forward Uncertainty Propagation}
Instead of relying on the sigma-point method required for exact Gaussian moment matching in (\ref{eq:general_moment_matching}), the analytical integration solutions from \cite{Deisenroth2015, Pan2018, Filabadi2025} for nonparametric GPs can be extended to semi-parametric models using basis functions with known Gaussian expectations and variances. Alternatively, as noted in Remark \ref{rem:sigma_points_forward_backward}, sigma points from the backward pass can be reused for forward propagation, avoiding extra computation.

\subsubsection{Online Sparsity Update}
We can consider two approaches for updating the online sparse semi-parametric GP. The first, proposed in \cite{Duy2011, He2024, Chowdhary2015, Csató2002}, uses a tolerance threshold on the elimination score to determine whether newly observed data should be removed, which may result in a variable-size data pool. 
In contrast, our elimination score depends on the parametric component, which captures information from all observations via the basis-function structure without increasing dimensionality. Consequently, it is not favourable to directly compare the elimination score with a fixed predefined threshold. Therefore, we update the GP after each observation and remove the data point with the lowest elimination score from the nonparametric data pool, thereby maintaining a fixed-size data pool equal in size to the initial dataset. 
This also aligns with modern just-in-time (JIT) compilation technology, which can leverage fixed shapes to improve computational efficiency.

\subsubsection{Hyperparameter Tuning}
The hyper-parameters of the kernel functions, $\kappa_i$, and the basis functions, $\phi_i$, can be tuned by maximizing the likelihood described in Remark \ref{rem:liklihood_theta}, using the initial dataset $\mathcal{D}_0$.
Specifically, the basis functions, $\phi_i$, can be any type of regressor, such as parametric physics-based models or neural networks, where the weights of the last layer can be fine-tuned during the online process, while the remaining parameters are optimized during hyperparameter tuning.

\section{Simulation Results} \label{sec:Results}

To evaluate the performance of the proposed methods, we conducted a series of experiments. 

We begin by describing the simulation settings, followed by experiments demonstrating the performance and capabilities of the stochastic FHDP introduced in Sec. \ref{sec:Traj_Opt} as an approximate method for stochastic optimal control. Subsequently, we present simulations applying the GPFHDP-AIF algorithm described in Sec. \ref{sec:algorithm}, first to a one-dimensional dynamical system and then to an autonomous vehicle system for motion planning. In the latter study, by incorporating AIF, the algorithm is able to steer the vehicle to concurrently generate and collect richer data for more efficiently learning the system dynamics, which leads to more accurate predictions for synthesizing optimal control policies in real time.  

All experiments were implemented in Python using JAX and executed on a 1.80 GHz Intel Core i7-1265U CPU. 
In all our experiments, we used a quadratic cost function of the following form over a horizon: 
\begin{subequations}\label{eq:quad_cost}
	\begin{align*}
		c_k(\xi_k) &= (x_k-x_k^{\text{ref}})^{\top} W (x_k-x_k^{\text{ref}}) + u_k^{\top} R u_k, \\
		c_{t+H}(x_{t+H}) &= (x_{t+H}-x_{t+H}^{\text{ref}})^{\top} W_{H} (x_{t+H}-x_{t+H}^{\text{ref}}), 
	\end{align*}
\end{subequations}
where $x_k^{\text{ref}}$ represents the goal or desired reference we aim to track and $W$, $W_H \in \mathbb{R}^{n_x \times n_x}$, and $R \in \mathbb{R}^{n_u \times n_u}$ represent weight matrices.
To generate the unit sigma points and their corresponding weights, we employed fifth-order unscented transforms (UT5) using the spherical cubature rule \cite{sarkka2023, Hassan2023}, which generates $N_{\xi} = 2n_{\xi}^2+1$ sigma points. To discretise the continuous dynamics, we used a fourth-order Runge–Kutta integration method and a zero-order hold for the control action.

\subsection{Stochastic Optimal Control for Pendulum Swing-Up}

This experiment evaluates the performance of the stochastic FHDP method described in Sec. \ref{sec:Traj_Opt}, in comparison with Stochastic DDP \cite{Theodorou2010_SDDP} and PDDD \cite{Pan2018, Pan2014}. PDDP has been compared with other prior works \cite{Pan2018}. 
We apply all methods to GP dynamic models trained on data from a single pendulum. Planning is performed over one horizon with $H = 50$, using a purely nonparametric GP and excluding the exploration cost function.
The task is to swing the surrogate GP models of the pendulum from the downward position $(\theta = 0)$ to the horizontal position $(\theta = \pi)$. The state is $x = (\theta, \dot{\theta})$ with $n_x = 2$ and $n_u = 1$. Model parameters are chosen from \cite{Filabadi2024, Hassan2023}. The sampling time is $0.1\unit{\second}$, with initial and goal states $x_0 = (0, 0)$ and $x^{\text{ref}} = (\pi, 0)$. We set $W = 10^{-1}\mathbb{I}$, $R = 10^{-1}$, $W_H = 10\mathbb{I}$, and $\Sigma_{xx,0} = 10^{-3}\mathbb{I}$.

\begin{figure}[!h]
  \centering
  \includegraphics[width=1\columnwidth]{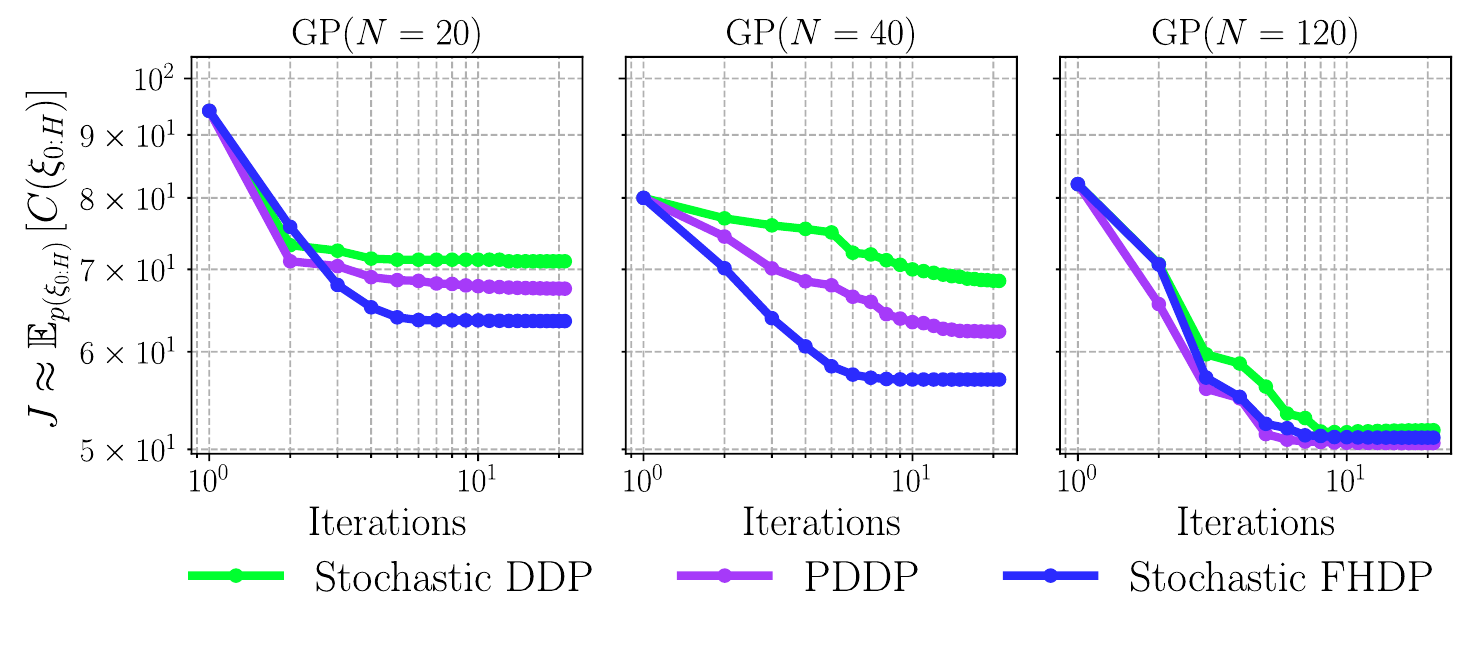} 
  \caption{\small{ 
        Convergence of the stochastic objective function computed over the normal distributions $\mathcal{N}(\xi_k; \mu_{\xi,k}, \Sigma_{\xi\xi,k})$ for $k=0,\dots,H$ over iterations for Stochastic DDP (green), PDDD (purple), and the proposed Stochastic FHDP (blue) using GP models trained with $N = 20, 40,$ and $120$ data points. 
  }}\label{fig:Stochastic_DDP_FHDP_PDDP}
\end{figure}

The evolution of the objective function (\ref{eq:SOC_objective}) over backward–forward iterations is shown in Fig. \ref{fig:Stochastic_DDP_FHDP_PDDP} for different GP model data pool sizes, $N\in\{20,40,120\}$.
Stochastic FHDP (blue) consistently achieves lower objective values and converges faster than both Stochastic DDP (green) and PDDD (purple).
This advantage is particularly significant for smaller data sets ($N=20$ and $N=40$), where the surrogate GP model is more stochastic (as noted in Lemma \ref{lem:GP_properties}) and the proposed method handles stochasticity more effectively.
As the number of data points increases ($N=120$), the stochasticity diminishes (in accordance with Lemma \ref{lem:GP_properties}), and the convergence behavior of all methods becomes similar. Also, the achieved ultimate optimal cost function correspondingly reduces from the left to the right subplot.

The average computation times are reported in Table \ref{tab:SOC_Time} for this experiment. 
Stochastic FHDP is consistently the fastest, followed by Stochastic DDP, while PDDP is slower and becomes much slower as the data size increases. 
The computational complexity of the different methods can be compared as follows.
As noted in \cite{Pan2018}, the forward pass of PDDP for a single-step exact moment matching has complexity $O(N^2 n_x^2 n_\xi)$, while the backward pass for single-step gradient evaluation is $O(N^2 n_b (n_b + n_u))$, where $n_b \coloneqq n_x + n_x(n_x+1)/2$ denotes the belief state dimension.
In contrast, the proposed Stochastic FHDP approach exhibits lower complexity. For the forward pass in (\ref{eq:general_moment_matching}), it is $O(N n_x N_\xi)$, and performing Cholesky factorization leads to $O(n_\xi^3)$. The backward pass, for computing the coefficients of the optimal state-action value function in (\ref{eq:FH_Q_1}) and (\ref{eq:FH_Q_3}), requires $O(N_\xi)$, assuming that the evaluations of the dynamics at sigma points are reused from the forward pass, as highlighted in Remark \ref{rem:sigma_points_forward_backward}.
For Stochastic DDP, following the analysis in \cite{Liao1991}, the computational complexity for gradient and Hessian evaluation with a GP dynamics model is $O(N(n_x^3 + n_x^2 n_u + n_x n_u^2))$. To ensure a fair comparison, we employed line search on the stochastic objective function while using the same uncertainty propagation as in (\ref{eq:general_moment_matching}) during the forward pass of Stochastic DDP, which introduces additional computational overhead.
The complexity of other algorithmic components is omitted, as they are similar for the different methods.

These results demonstrate that Stochastic FHDP provides a more reliable and efficient trajectory optimization approach for stochastic settings than the baseline methods.

\begin{table}[!h] 
\centering
\begin{tabular}{c||ccc}
\hline \hline
\multirow{2}{*}{\makecell{No. Data \\ in Data Pool}}  & \multicolumn{3}{c}{Average run-time per Iteration [ms]} \\ \cline{2-4}
& \makecell{Stochastic FHDP} 
& \makecell{Stochastic DDP} 
& \makecell{PDDP} 
\\ \hline  \hline
\multirow{1}{*}{$N=20$} 
& 4.8 & 5.4 & 28.6  \\ \hline
\multirow{1}{*}{$N=40$} 
& 7.1 & 11.5 & 78.9   \\ \hline
\multirow{1}{*}{$N=120$} 
& 40.9 & 110.2 & 863.1   \\ \hline
\hline
\end{tabular}
\caption{
\small{ The average run-time for different algorithms.}}\label{tab:SOC_Time}
\end{table}

\subsection{GPFHDP-AIF for 1D case}
\begin{figure}[!h]
  \centering
  \includegraphics[width=1\columnwidth]{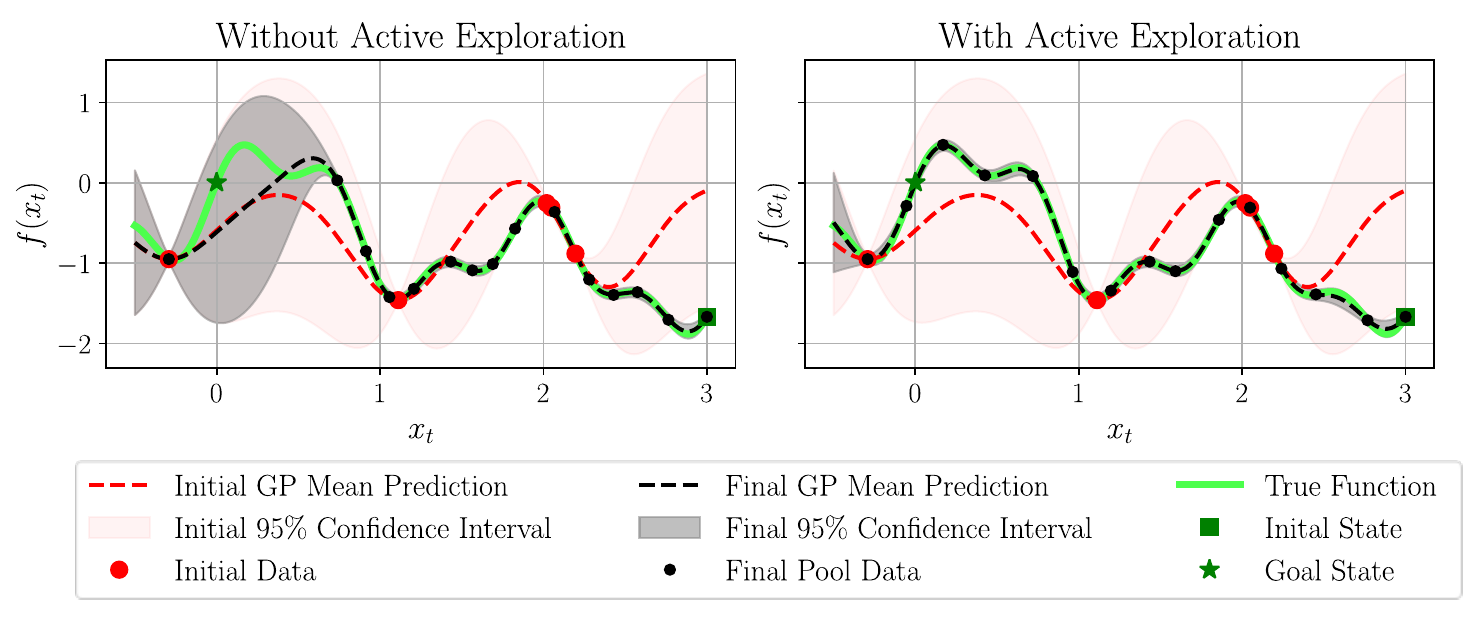} 
  \caption{\small{True dynamic function model (green) compared with the initial (red) and final surrogate model (black) after 40 time steps, with and without active exploration. The initial dataset $\mathcal{D}_0$ is shown as red dots, and the final stored data as black dots. The system starts from the initial state (green square) and moves towards the goal state (green star) using the proposed algorithm.  
  }}\label{fig:ActExp_1DCase}
\end{figure}

To illustrate the properties of the proposed algorithm, we applied it to the following one-dimensional system and visualize the agent’s behaviour over the function domain in Fig. \ref{fig:ActExp_1DCase} with and without using the exploration cost term. 

\begin{figure}[!h]
  \centering
  \includegraphics[width=0.75\columnwidth]{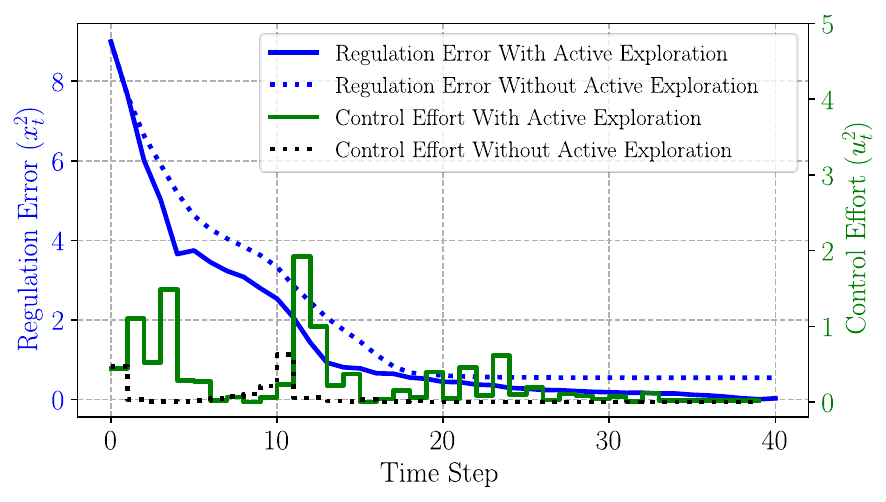} 
  \caption{\small{ Regulation error and control effort for the one-dimensional dynamics, with and without active exploration. 
  }}\label{fig:ActExp_1DCase_XU}
\end{figure}

In this experiment, we used a sparse online GP without a parametric component, starting with 5 data points in the GP data pool and storing a maximum of 15 data points. Therefore, we incrementally increase the data pool size up to 15 and then keep it fixed, rather than maintaining a fixed size throughout the experiment.

\begin{multline*}
    f(x) = \tanh\left(1 + 0.05x - 0.5x^{2}\right) + 0.6  \sin(4x) \\
           \quad + 0.3  \sin(10x + 0.5) e^{-0.05x} -0.14 
\end{multline*}
\begin{equation*}
    \dot{x}_t = f(x_t) + u_t 
\end{equation*}
This experiment aims to move the system from the initial state $x_0 = 3$ to the origin as the target state, $x^{\text{ref}} = 0$. The regulation error and control effort over time are shown in Fig. \ref{fig:ActExp_1DCase_XU}.
As shown in Fig. \ref{fig:ActExp_1DCase}, incorporating the exploration cost significantly reduces both the uncertainty and the error of the final surrogate model after 40 time steps, compared to the case without active exploration. Without active exploration, the agent fails to reach the goal state, and the uncertainty around the goal remains high. This effect is also evident in Fig. \ref{fig:ActExp_1DCase_XU}, where the regulation error converges to zero under active exploration, with a faster convergence rate, although requiring greater excitation in the control input.

\subsection{Autonomous Vehicle Experiment}

To evaluate our algorithms, we consider a 6-DoF model to represent vehicle dynamics \cite{JeanPierre2022, Chuan2021}. The states are global position coordinates $X, Y$, yaw angle $\psi$, longitudinal and lateral velocities $v_x, v_y$, and yaw rate $\omega$. The control inputs are the front wheel steering angle $\delta$ and normalized longitudinal tire force at the front wheel $\tilde{F}_{x,f}$. The body reference frame is centered at the vehicle’s center of gravity, with the X-axis oriented toward the front. The kinematic equations for $X$, $Y$, and $\psi$ are expressed in the global Cartesian frame, while forces and moments are modeled in the body frame. 
The full vehicle dynamics, slip angle formulations, and model parameters are given in \cite{Filabadi2025}. The state and control vectors are defined as $x \coloneqq \begin{bmatrix} X & Y & \psi & v_x & v_y & \omega\end{bmatrix}^{\top}$ and $u \coloneqq \begin{bmatrix} \delta & F_{x,f} \end{bmatrix}^{\top}$. With $n_x = 6$, $n_u = 2$, and $n_{\xi} = 8$.
The sampling time is set to $0.1\unit{\second}$ following \cite{Askari2025}.
Next, we begin with experiments on the proposed sparse online semi-parametric model learning, then we demonstrate the results obtained with the GPFHDP-AIF algorithm on the simulated vehicle.

\subsubsection{Sparsification Experiment}
In this section, we evaluate the impact of using different types of basis functions within the sparsification method proposed in Sec. \ref{sec:Inf4Model}. We also examine the role of the parametric component alongside the nonparametric component.
The basis functions considered are detailed in Table \ref{tab:basis_fn}. For the experiments, we set $d=5$, $8$, and $16$ for the Polynomial base, Fourier base, and Radial Basis Function (RBF), respectively.

\begin{table}[!h] 
\centering
\begin{tabular}{c||ccc}
\hline \hline
\makecell{Basis \\ Function} & $\phi_i(\xi)$ & \makecell{No. \\ Units} & \makecell{No. Hyper-\\ Parameters} 
\\ \hline  \hline
Linear
& $\tanh[1, \xi]$ & $n_{\xi}+1$ & -  \\ \hline
Polynomial 
& $\tanh[1, \xi, \xi^2, \dots ,\xi^d]$ & $dn_{\xi}+1$ & -   \\ \hline
Fourier 
& $[1, \dots, \sin(w_d\xi), \cos(w_d\xi)]$ & $2d+1$ & $d$  \\ \hline
RBF 
& $\left[1, \dots, \exp\left(- \tfrac{\lVert \xi-c_d \rVert^2}{l_d} \right)\right]$ & $d+1$ & $d+n_{\xi}$   \\ \hline
\hline
\end{tabular}
\caption{
\small{ Basis Functions.}}\label{tab:basis_fn}
\end{table}

\begin{figure}[!h]
    \centering
    \begin{subfigure}[b]{0.5\textwidth}
        \centering
        \includegraphics[width=0.75\linewidth]{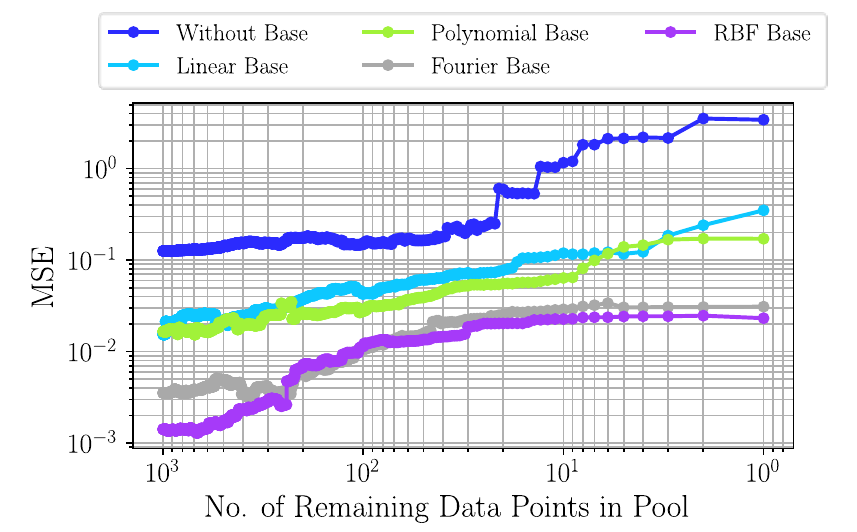}
        \caption{}
        \label{fig:Sparse_GP_MSE}
    \end{subfigure}
    \hfill
    \begin{subfigure}[b]{0.5\textwidth}
        \centering
        \includegraphics[width=0.75\linewidth]{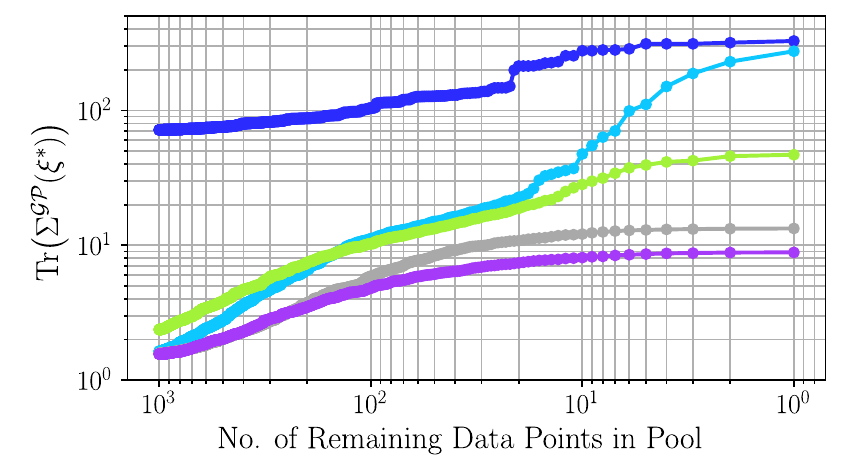}
        \caption{}
        \label{fig:Sparse_GP_Cov}
    \end{subfigure}
    \caption{\small{
    MSE (a) and predictive uncertainty (b) of different semi-parametric and nonparametric models as the number of remaining data points in data pool decreases.}
    }
    \label{fig:Sparse_GP}
\end{figure}

We randomly selected 1000 data points for training the GP and 64 data points as the test set, drawn from the same distribution as the training data. 
The Adam optimizer was employed to tune the hyperparameters $\phi_i$ and $\kappa_i$. After initializing the parameters in (\ref{eq:postrior_params_GP}) and (\ref{eq:def_recursive_param}), we applied the sparsification method described in Lemma \ref{lem:Sparse_SemiParametric_GP} to sequentially remove data points from the pool used by the nonparametric component. At each step, we computed the predictive mean and covariance of the resulting GP on the test set. Finally, we evaluated the prediction performance of the semi-parametric GP by demonstrating the MSE and the trace of the predictive covariance in Fig. \ref{fig:Sparse_GP}. In this experiment, the parametric part of the models was trained on the entire training set, so the effects of removing data from the nonparametric part were investigated, along with different basis functions.

As shown in Figs. \ref{fig:Sparse_GP_MSE} and \ref{fig:Sparse_GP_Cov}, all methods exhibit increasing error and uncertainty as fewer data points remain, consistent with the continuity properties of the GP covariance described in Lemma \ref{lem:GP_properties}. Semi-parametric models equipped with basis functions consistently outperform the nonparametric model without any base. In particular, the RBF basis achieves the lowest MSE and uncertainty across all models, followed by the Fourier basis, while linear and polynomial bases provide moderate improvements.
According to Fig. \ref{fig:Sparse_GP}, semi-parametric models with fewer than 10 data points perform equivalently to a purely nonparametric model trained on 1000 data points. This demonstrates that incorporating a parametric component can drastically reduce computational effort by requiring fewer data points. Moreover, basis functions with tunable hyperparameters, the RBF and Fourier bases, capture more information and offer the most robust performance under data sparsity.

\subsubsection{GPFHDP-AIF Experiment}

In this section, we apply the GPFHDP-AIF algorithm to vehicle dynamics. We use the same RBF parametric model described in Table \ref{tab:basis_fn}, whose performance was evaluated in the previous section. For the receding horizon implementation, we set $H = 10$ and $\Sigma_{xx,0} = 10^{-3} \mathbb{I}$. The vehicle starts from an initial longitudinal velocity $V_0$ and with the initial state $x_0 =  \begin{bmatrix} 0 & 0 & 0 & V_0 & 0 & 0 \end{bmatrix}^{\top}$. 
To evaluate vehicle manipulation under the AIF framework, we consider two scenarios: a regulation scenario, representing lane-changing or overtaking manoeuvres, and a tracking scenario, where the vehicle follows a predefined infinity-shaped path.
We investigate the performance of the semi-parametric GP using different data pool sizes, $N \in \{10, 30, 60\}$, particularly when incorporating the active exploration cost with a regularization parameter $\gamma = 5$. For benchmarking, we compare the results against an optimal model-based controller, implemented as a receding horizon version of deterministic FHDP with sufficient iterations at each time step.
In all experiments, the initial dataset $\mathcal{D}_0$ is generated randomly, with its size equal to that of the data pool, and the hyperparameters are tuned using this set.

\begin{figure}[!h]
  \centering
  \includegraphics[width=1.0 \columnwidth]{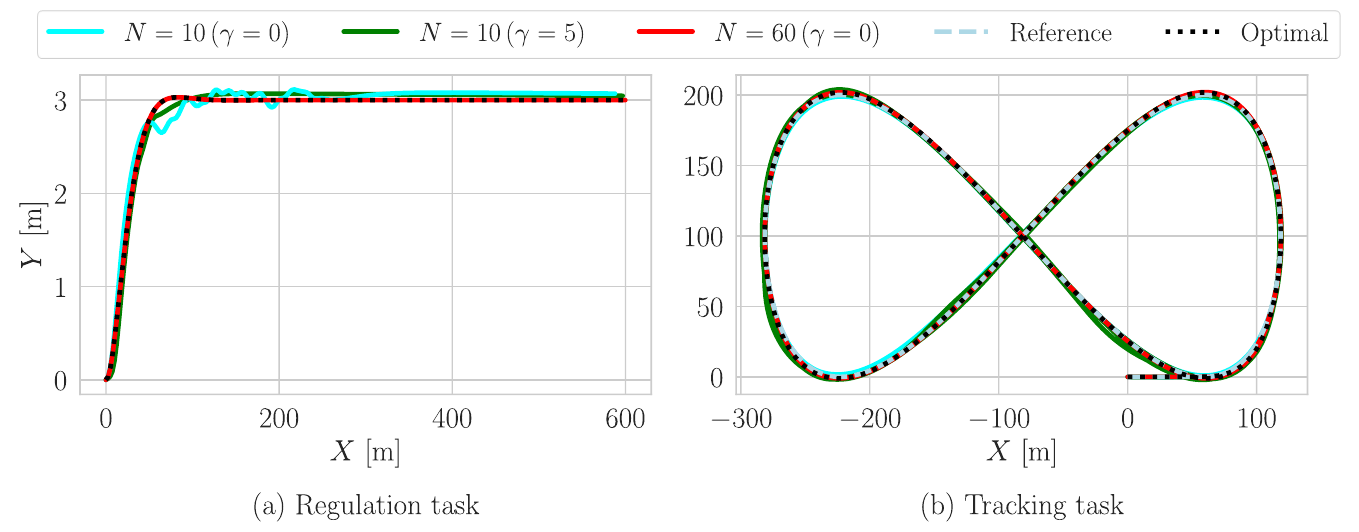} 
  \caption{\small{ X–Y trajectories for $N = 10$, $N = 60$, and the optimal trajectory for (a) the regulation task and (b) the tracking task.
  }}\label{fig:ActExp_Car_XY}
\end{figure}

\begin{figure}[!h]
  \centering
  \includegraphics[width=1.0\columnwidth]{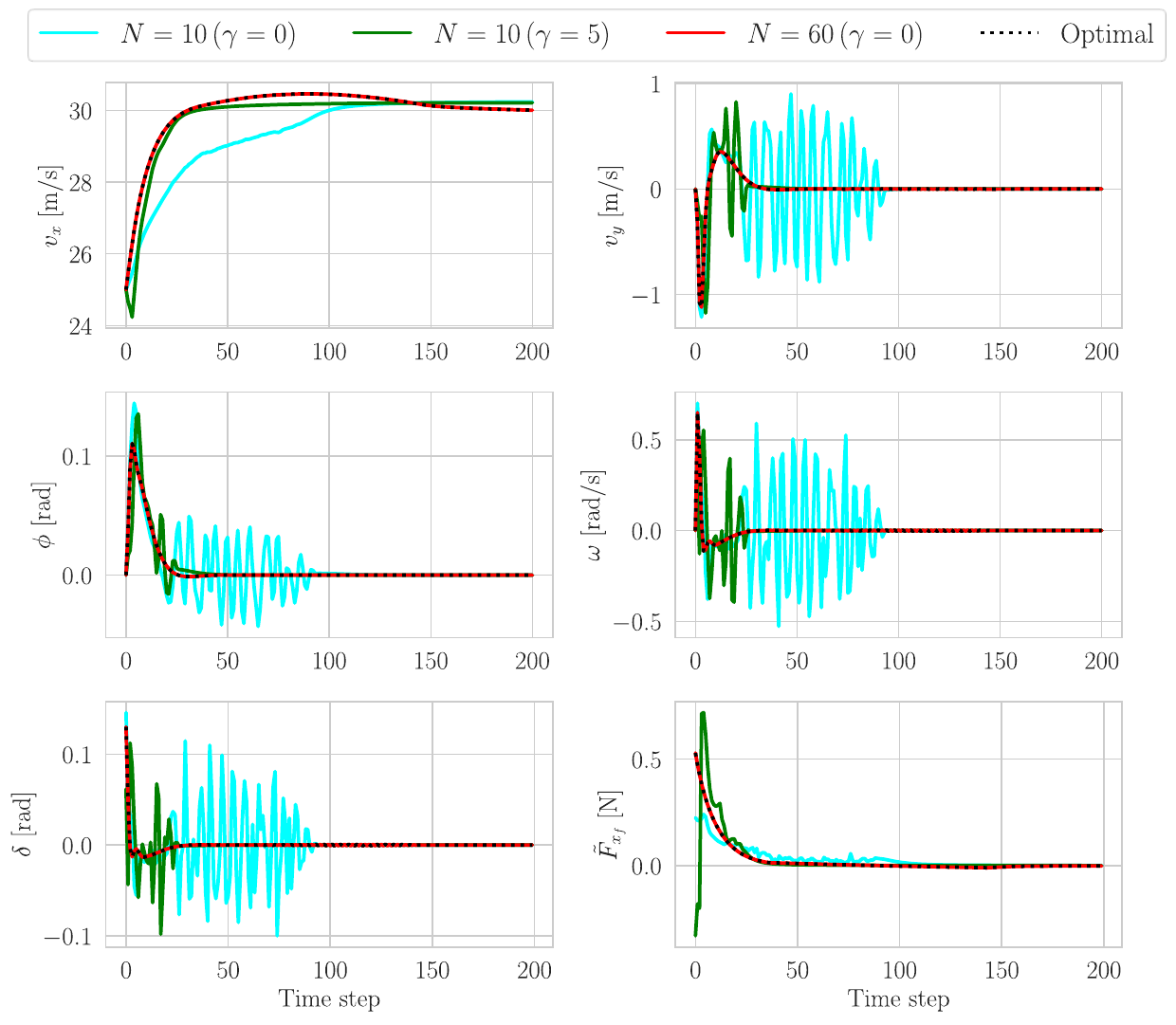} 
  \caption{\small{State and control trajectories of the optimal model-based algorithm and GPFHDP-AIF algorithm for the regulation task. 
  }}\label{fig:ActExp_Car_Reg_State_Control}
\end{figure}

\begin{figure}[!h]
  \centering
  \includegraphics[width=1.0\columnwidth]{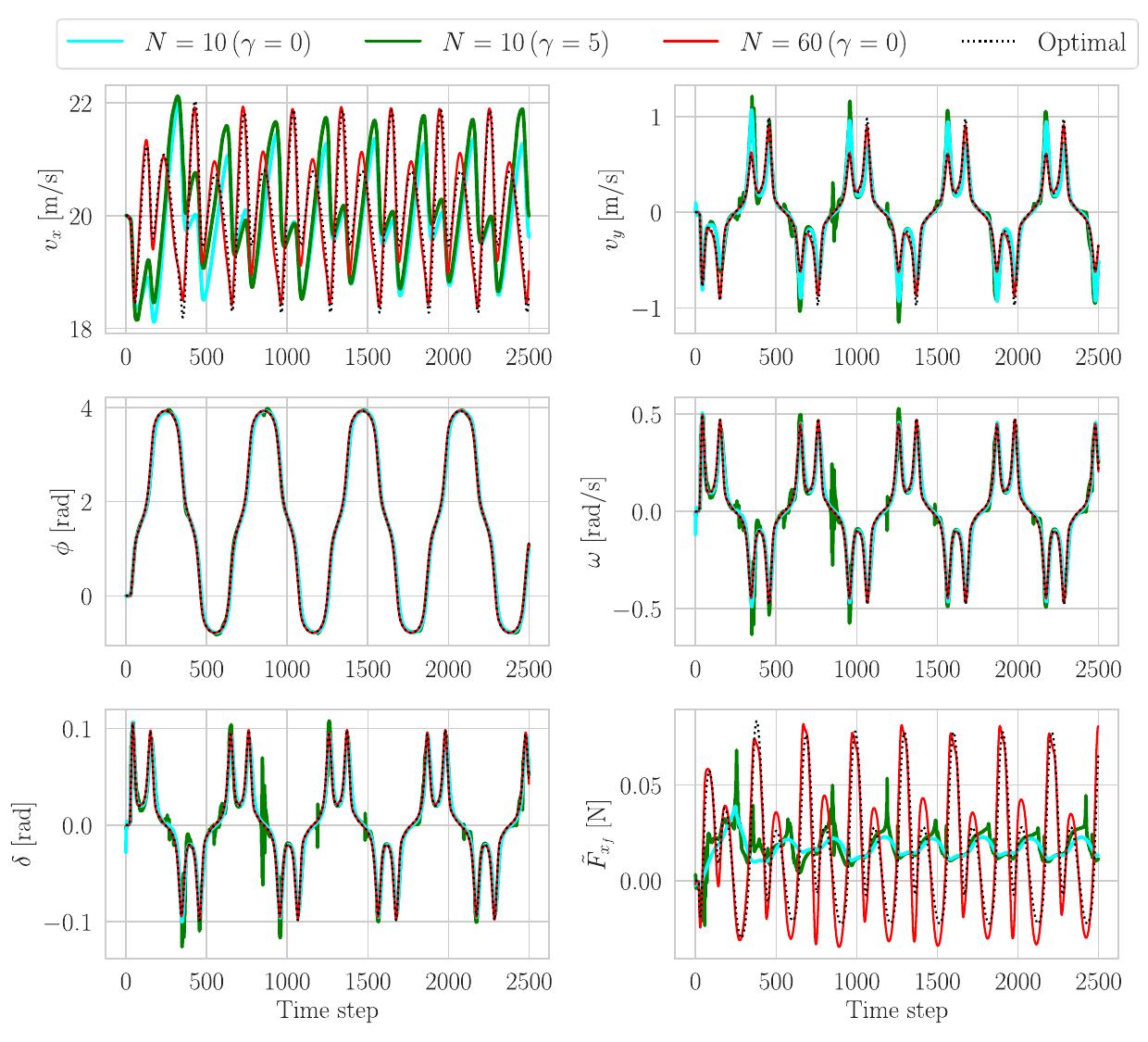} 
  \caption{\small{State and control trajectories of the optimal model-based algorithm and GPFHDP-AIF algorithm for the tracking task. 
  }}\label{fig:ActExp_Car_Tracking_State_Control}
\end{figure}

The results for all scenarios are summarized in Table \ref{tab:results_vehicle}. Trajectories on the X–Y plane are shown in Fig. \ref{fig:ActExp_Car_XY}, while the state and control trajectories are presented in Fig. \ref{fig:ActExp_Car_Reg_State_Control} and Fig. \ref{fig:ActExp_Car_Tracking_State_Control} for the regulation and tracking tasks, respectively.
As presented in Table \ref{tab:results_vehicle}, increasing the data pool size progressively yields results that approach the optimal solution, although at the cost of longer computation times. For instance, with $N = 60$, a near-optimal solution can be achieved. Notably, the computation time consistently remains below the sampling time, thereby ensuring the feasibility of real-time processing.
Using the active exploration term ($\gamma = 5$) instead of removing it ($\gamma = 0$) generally reduces the error and gives more optimal solutions. This effect becomes especially clear when the data pool is small or, equivalently, when the mismatch between system and surrogate model is large, such as in the case of $N = 10$.

\begin{table*}[h!]
\centering
\resizebox{0.7\textwidth}{!}{
\begin{tabular}{c||c|ccc||cccc}
\hline \hline
\multirow{3}{*}{\makecell{No. Data \\ in Data Pool}}  & \multirow{3}{*}{$\gamma$} & \multicolumn{3}{c||}{Regulation Task} & \multicolumn{3}{c}{Tracking Task} \\ \cline{3-8}
&&\makecell{Error} 
& \makecell{Control \\ Effort} 
& \makecell{Average run-time \\ per time-step [ms]} 
& \makecell{Error} 
& \makecell{Control \\ Effort} 
& \makecell{Average run-time \\ per time-step [ms]} 
\\ \hline  \hline
\multirow{2}{*}{$N=10$} 
& 0 & 65.36 & 3.68 & 18.03 & 91.33 & 1.09 & 21.03 \\ \cline{2-8}
& 5 & 62.68 & 5.08 & 18.62 & 72.55 & 1.08 & 20.02 \\  \hline
\multirow{2}{*}{$N=30$} 
& 0 & 63.37 & 6.13 & 46.09 & 48.72 & 1.22 & 51.02 \\ \cline{2-8}
& 5 & 61.93 & 3.57 & 45.62 & 41.78 & 1.24 & 52.02 \\  \hline
\multirow{1}{*}{$N=60$} 
& 0 & 61.81 & 3.82 & 95.92 & 5.83 & 1.22 & 98.92 \\ \hline
\makecell{Optimal Model-Based}
& -- & 61.53 & 3.80 & 9.51 & 5.78 & 1.20 & 10.01 \\ 
\hline \hline
\end{tabular}
}
\caption{\small{Results of the GPFHDP-AIF algorithm on vehicle dynamics.}}\label{tab:results_vehicle}
\end{table*}

To evaluate performance, we compare the trajectories obtained with $N = 10$ and $N = 60$ with the optimal trajectories, as shown in Figs. \ref{fig:ActExp_Car_Reg_State_Control} and \ref{fig:ActExp_Car_Tracking_State_Control}.
In Fig. \ref{fig:ActExp_Car_Reg_State_Control}, which corresponds to the regulation task, the agent without active exploration explores the dynamics in an unguided manner to reach the goal. In contrast, with active exploration the agent engages in goal-oriented exploration, learns from more informative data, and therefore converges to the goal state more quickly.

In Fig. \ref{fig:ActExp_Car_Tracking_State_Control}, which corresponds to the tracking task, it can be seen that without using the active exploration term, the agent simply attempts to follow the path. By contrast, incorporating active exploration produces initial excitations during the early time steps. After passing the path for some rounds, model exploitation begins, the oscillations diminish, and the agent approaches a sub-optimal solution.

\begin{figure}[!h]
  \centering
  \includegraphics[width=0.7\columnwidth]{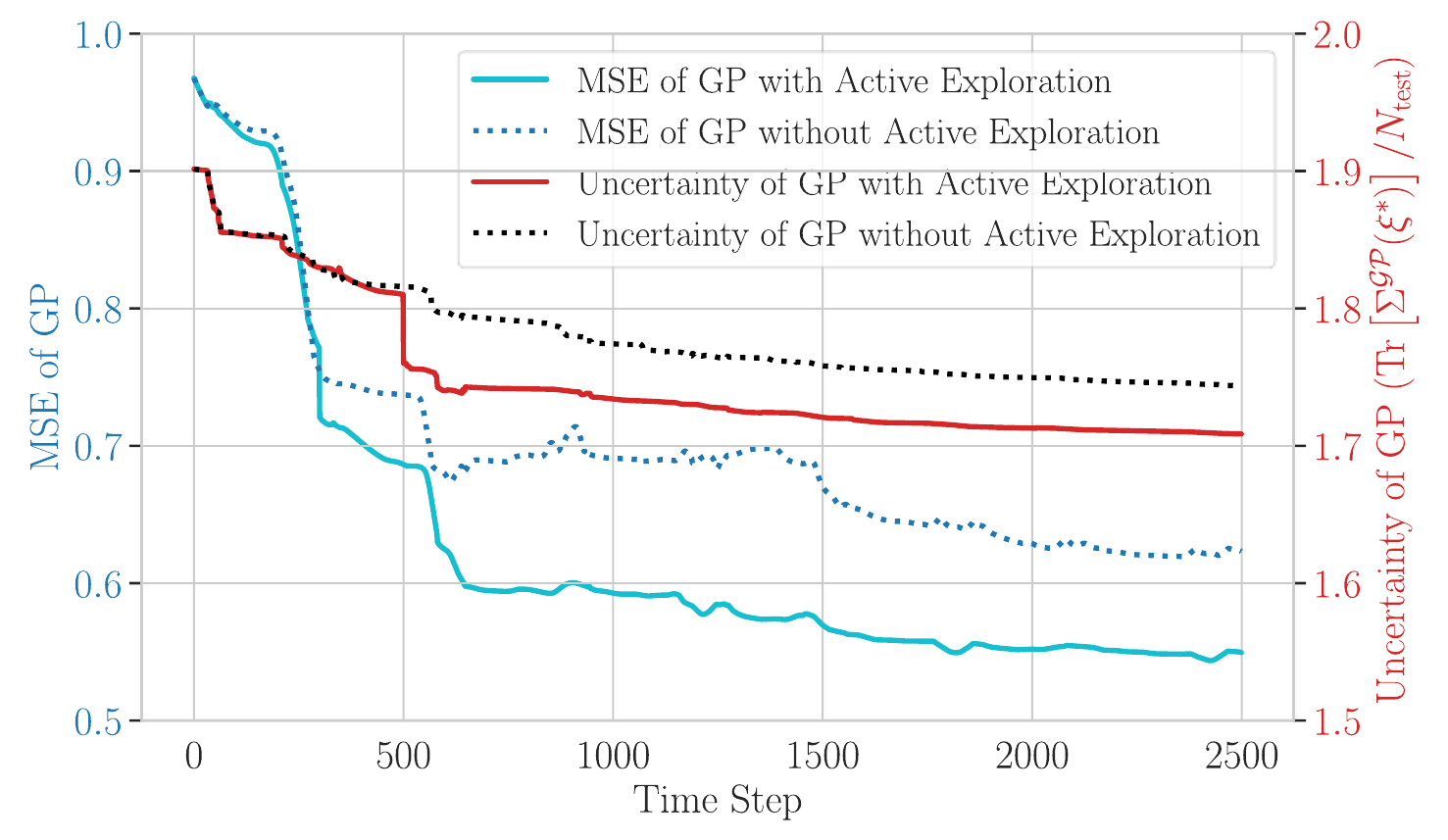} 
  \caption{\small{Mean squared error and uncertainty of the GP surrogate model evaluated on a test set at each time step for the tracking task, with and without active exploration. 
  }}\label{fig:ActExp_Tracking_MSE_UncertaintyGP}
\end{figure}

\begin{figure}[!h]
  \centering
  \includegraphics[width=0.7\columnwidth]{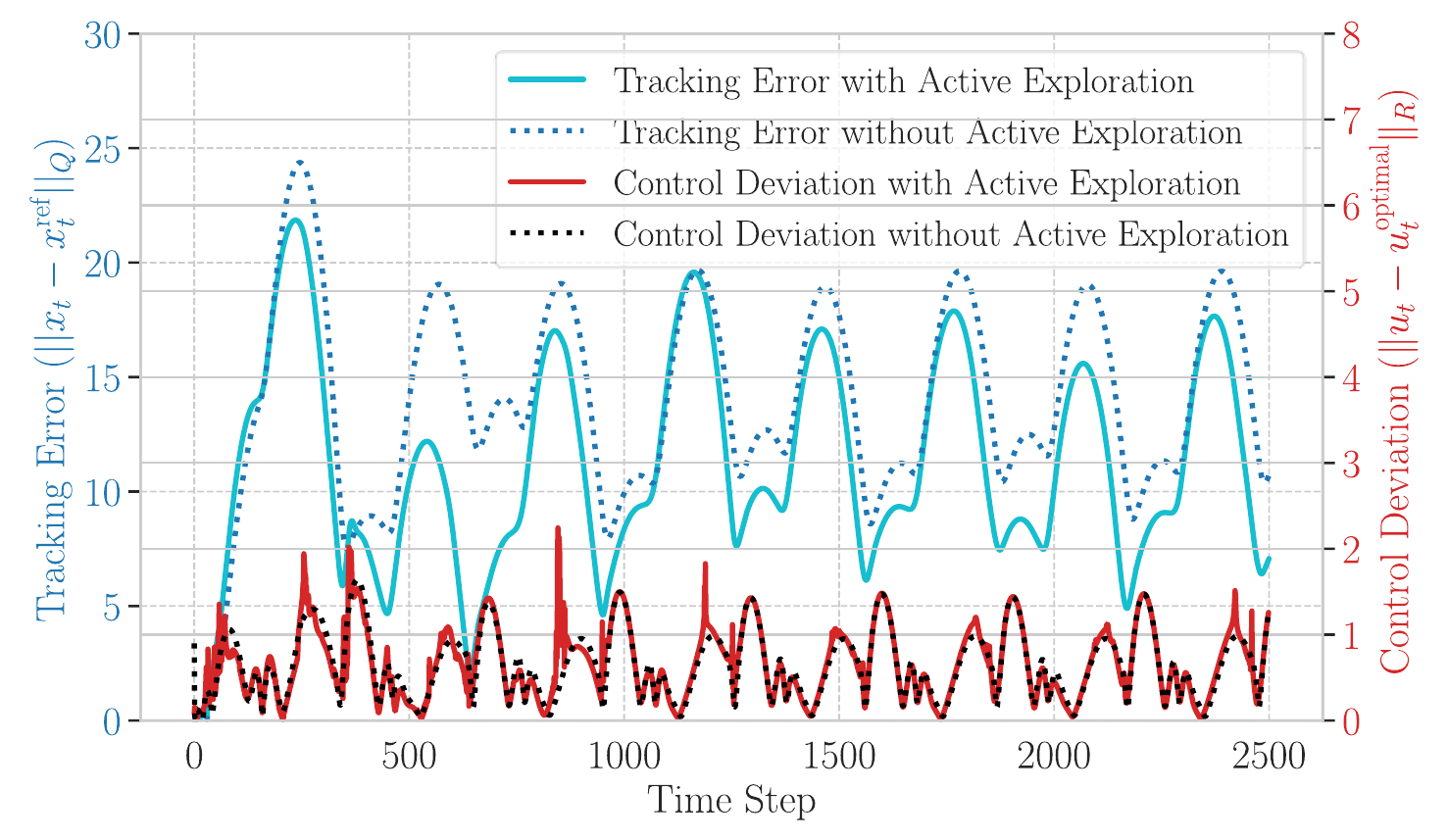} 
  \caption{\small{Deviation of position and control input from the reference trajectory and the optimal solution, with and without active exploration for the tracking task.
}}\label{fig:ActExp_Tracking_TrackingError_ControlEffort}
\end{figure}

Additional details for this scenario are provided in Figs. \ref{fig:ActExp_Tracking_MSE_UncertaintyGP} and \ref{fig:ActExp_Tracking_TrackingError_ControlEffort}. Fig. \ref{fig:ActExp_Tracking_MSE_UncertaintyGP} presents the evaluation of model learning at each time step on a test dataset of 1000 samples. These results indicate that active exploration helps reduce model error and uncertainty largely. In other words, as more informative data is collected and the model is updated, uncertainty decreases further, enabling the algorithm to converge to a better solution more reliably.

Fig. \ref{fig:ActExp_Tracking_TrackingError_ControlEffort} shows that due to exploration, the agent exhibits some oscillations in control input during the initial time steps but achieves a lower tracking error. After a number of time steps, the control deviation becomes comparable to the case without active exploration, while the tracking error remains reduced. In other words, the method invests more control effort in the short term to achieve improved long-term tracking performance.

Overall, in the regulation task, the greater freedom for exploration enables the system to reach the optimal solution within only a few time steps, whereas in the tracking task several cycles are needed to approach a more optimal solution.

\section{Conclusion} \label{sec:Conclusion}

In this work, we presented a principled and computationally tractable solution for dual MPC that integrates a real-time probabilistic model learner with AIF to tackle the fundamental challenge of balancing information gathering (identification and exploration) with maintaining performance (optimization and exploitation). Our approach explicitly considers the stochastic nature of the resulting decision-making problem in the dual control design. 

The proposed framework led to an SOC problem defined over a receding horizon in the presence of a novel online-updating probabilistic model. The cost function integrates control performance with an additive exploration term, encouraging the controller to generate informative data while maintaining overall performance. To solve the resulting SOC, we employed an iterative method based on local approximations of dynamic programming.
Simulation experiments demonstrated the capability of the proposed algorithm to design intelligent controllers, particularly for autonomous vehicles.   

Future research directions include extending the framework to partially observable systems and incorporating safety constraints into the proposed algorithm.

\appendix  

\subsection{Proof of Lemma \ref{lem:GP_properties}} \label{app:lem1}
The posterior covariance of the GP defined in (\ref{eq:GP_semiparametric}) and mentioned in (\ref{eq:pred_SemiGP}) equivalently can be expressed as:
    \begin{equation*} 
        \Sigma^{\mathcal{GP}}_{ii,t} (\xi^*) = \mathcal{\tilde{K}}_{i}^{**}(\xi^*) - \mathcal{\tilde{K}}_{i,t}^{*,\top}(\xi^*) \mathcal{\tilde{K}}_{i,t}^{-1} \mathcal{\tilde{K}}_{i,t}^{*}(\xi^*).
    \end{equation*} 
    The second term is always non-negative, and a trivial bound is given by $\mathcal{\tilde{K}}_{i}^{**}(\xi^*)$. One can obtain 
    \begin{equation*} 
        \Sigma^{\mathcal{GP}}_{ii,t} (\xi^*) \le \tilde{\kappa}_i (\xi^*,\xi^*) \le A_i + \sigma^2_{f_i}+ \lambda_{\text{max}}[\Sigma_{\theta_i,0}] \bar{\phi}^2_{i}.  
    \end{equation*} 
    This bound reflects the fact that conditioning on data always reduces variance \cite{Rasmussen2006}. 
    Additionally, if Assumption \ref{assum:RKHS} holds, then by applying the consistency property mentioned in \cite{Rasmussen2006} on the GP defined in (\ref{eq:GP_semiparametric}), (\ref{eq:GP_consistency}) will hold, as also noted in \cite{Alpcan2015}. 

\subsection{Proof of Lemma \ref{lem:Online_SemiParametric_GP}} \label{app:lem2}
When a new data pair $(\xi_t, x_{t+1})$ becomes available, we can partition the matrices $Y_{i,t+1}$, $\Phi_{i,t+1}$, and $\mathcal{K}_{i,t+1}$ as follows:
    \begin{equation}
            \mathcal{K}_{i,t+1} = \begin{bmatrix}              				
                \mathcal{K}_{i,t} & \mathcal{K}^{*}_{i,t} (\xi_{t})  \\
                 		    \square & \mathcal{K}^{**}_{i,t} (\xi_{t}) \\
                    			    \end{bmatrix}, 
    \end{equation}
    \begin{equation}\label{eq:Y_Phi_partitioned}
            Y_{i,t+1} = [Y^{\top}_{i,t}, x_{i,t+1}]^{\top}, \quad
            \Phi_{i,t+1} = [\Phi_{i,t}, \Phi^{*}_{i,t}(\xi_{t})].
    \end{equation}
    Applying the block matrix inversion lemma for $\mathcal{K}_{i,t+1}$, one can obtain 
    \begin{equation} \label{eq:K_inv}
                \begin{aligned}
                \mathcal{C}_{i,t+1} &= \mathcal{K}_{i,t+1}^{-1} \\
                &= \begin{bmatrix}              				
                 \mathcal{K}_{i,t}^{-1} & 0  \\
                 		         0 & 0 \\
            \end{bmatrix}
            + b_{i,t}
            \begin{bmatrix}              				
                \mathcal{K}_{i,t}^{-1}\mathcal{K}^{*}_{i,t} (\xi_{t})   \\
                 		    -1  \\
            \end{bmatrix}
            \begin{bmatrix}              				
                \mathcal{K}_{i,t}^{-1}\mathcal{K}^{*}_{i,t} (\xi_{t})   \\
                 		    -1  \\
            \end{bmatrix}^{\top}\\
                &= \mathcal{U} [\mathcal{C}_{i,t}] + b_{i,t} S_{i,t} S^{\top}_{i,t},
                \end{aligned}
    \end{equation}
    which yields the recursion in (\ref{eq:GP_recursive_update_C}). By constructing $\alpha_{i,t+1}$ and $\beta_{i,t+1}$ according to the definitions in (\ref{eq:def_recursive_param}), using the partitioned matrices in (\ref{eq:Y_Phi_partitioned}) and the decomposed matrix $\mathcal{K}_{i,t+1}^{-1}$ from (\ref{eq:K_inv}), we can obtain the recursions presented in (\ref{eq:GP_recursive_update_alpha}) and (\ref{eq:GP_recursive_update_beta}). Moreover, we can expand (\ref{eq:postrior_params_GP_sigma}) to obtain a real-time update as
    \begin{align*}
        \Sigma_{\theta_i,t+1} &= \left(\Phi_{i,t+1}\mathcal{K}_{i,t+1}^{-1} \Phi_{i,t+1}^{\top} + \Sigma_{\theta_i,0}^{-1}\right)^{-1} \\
        &= \left(\Phi_{i,t}\mathcal{K}_{i,t}^{-1} \Phi_{i,t}^{\top} + b_{i,t} \mathcal{R}_{i,t}\mathcal{R}_{i,t}^{\top} + \Sigma_{\theta_i,0}^{-1}\right)^{-1}\\
        &= \left( \Sigma_{\theta_i,t}^{-1} + b_{i,t} \mathcal{R}_{i,t}\mathcal{R}_{i,t}^{\top}\right)^{-1} \\
        &= \left(\mathbb{I} - \Sigma_{\theta_i,t} \mathcal{R}_{i,t} (1/b_{i,t}+\mathcal{R}_{i,t}^{\top}\Sigma_{\theta_i,t}\mathcal{R}_{i,t}) ^{
        -1
        }\mathcal{R}_{i,t}^{\top}\right)\Sigma_{\theta_i,t} \\
        &= \left(\mathbb{I}_{n_{\theta_i}} - \Lambda_{i,t}\mathcal{R}_{i,t}^{\top} \right) \Sigma_{\theta_i,t},
    \end{align*}
    where the second equality follows from (\ref{eq:Y_Phi_partitioned}) and (\ref{eq:K_inv}), the third from (\ref{eq:postrior_params_GP_sigma}), and the fourth from applying the Woodbury matrix identity. Here, $\Lambda_{i,t}$ is defined in (\ref{eq:GP_recursive_update_Lambda}).
    Similarly, for (\ref{eq:postrior_params_GP_mu}), we obtain the following real-time update rule:
    \begin{align*}
        \mu_{\theta_i,t+1} &= \Sigma_{\theta_i,t+1}\left(\Phi_{i,t+1} \mathcal{K}_{i,t+1}^{-1}Y_{i,t+1} + \Sigma_{\theta_i,0}^{-1}\mu_{\theta_i,0}\right) \\
        &= \Sigma_{\theta_i,t+1}\left(\Phi_{i,t} \mathcal{K}_{i,t}^{-1}Y_{i,t} +a_{i,t}\mathcal{R}_{i,t} + \Sigma_{\theta_i,0}^{-1}\mu_{\theta_i,0}\right) \\
        &= \Sigma_{\theta_i,t+1}\left(\Sigma_{\theta_i,t}^{-1}\mu_{\theta_i,t} +a_{i,t}\mathcal{R}_{i,t} \right) \\
        &= \mu_{\theta_i,t} - \Lambda_{i,t}\mathcal{R}^{\top}_{i,t} \mu_{\theta_i,t} + a_{i,t} \Sigma_{\theta_i,t+1} \mathcal{R}_{i,t}\\
        &= \mu_{\theta_i,t} + \Lambda_{i,t} \left(-\mathcal{R}^{\top}_{i,t} \mu_{\theta_i,t} + a_{i,t} / b_{i,t} \right)\\
        &= \mu_{\theta_i,t} + \Lambda_{i,t}\left(\mu^{\mathcal{GP}}_{i,t} - x_{i,t+1}\right).   
    \end{align*}

\subsection{Proof of Lemma \ref{lem:Sparse_SemiParametric_GP}} \label{app:lem3}
Adopting the approach from \cite{Csató2002, Chowdhary2015}, first note that if the position of a data point in the dataset $\mathcal{D}^i_{t+1}$ is changed, the corresponding elements in the vectors $\alpha_{i,t+1}$ and $\beta_{i,t+1}$ will be changed accordingly. The same rule applies to the matrix $\mathcal{K}_{i,t+1}$, where reordering a data point’s position is represented by a similarity transformation $P^{\top} \mathcal{K}_{i,t+1} P$ in which $P$ is a permutation matrix that reorders specific rows and columns while leaving the remaining entries unchanged. Because $P$ is orthogonal and satisfies $P^{\top} = P^{-1}$, this reordering property also holds for the matrix $\mathcal{C}_{i,t+1} = \mathcal{K}^{-1}_{i,t+1}$. 
    Therefore, we start by proving the expressions in (\ref{eq:elimination_score}) and (\ref{eq:deletion_update}) for the last data point $(\xi_t, x_{i,t+1})$ appended to the dataset, i.e. $j = t+1$. Following the rule for changing a data point’s position, any elements corresponding to its data point in the vectors and matrices $\alpha_{i,t+1}$, $\beta_{i,t+1}$, $\mathcal{K}_{i,t+1}$, and $\mathcal{C}_{i,t+1}$ can be moved to the last row and column, and the same expressions remain valid. 

    We then show the computation of the elimination score in (\ref{eq:elimination_score}) for the last data point $(\xi_t, x_{i,t+1})$ appended to the dataset. The elimination score is defined as the RKHS norm of the difference between the mean functions of the GP posterior in (\ref{eq:GP_reparameterization}), with and without incorporating the newly added data point $(\xi_t, x_{i,t+1})$. It is important to note, based on Remark \ref{rem:GP_recursions}, while the new data point is included when updating $\mu_{\theta_i,t+1}$, it may still be eliminated from the dataset for the nonparametric component.
    \begin{equation} \label{eq:elimination_score_proof}
        \begin{multlined}
            \varepsilon_{i,t+1}(t+1) \\
            \resizebox{1\linewidth}{!}{$
            \begin{aligned}
                &\coloneqq \lVert \mu^{\mathcal{GP}}_{i,t+1}(\xi^*) - \mu^{\mathcal{GP}}_{i,t}(\xi^*)|_{\mu_{\theta_i,t+1}} \rVert_{\mathcal{H}}\\  
                &\begin{multlined}
                    = \lVert(\alpha_{i,t+1}-\beta_{i,t+1}\mu_{\theta_i,t+1})^{\top}\mathcal{K}_{i,t+1}^{*}(\xi^*) \\
                    - (\alpha_{i,t} - \beta_{i,t}\mu_{\theta_i,t+1})^{\top}\mathcal{K}_{i,t}^{*}(\xi^*) \rVert_{\mathcal{H}}
                    \end{multlined}\\
                 &\begin{multlined}
                    = \lVert (\mathcal{T} [\alpha_{i,t}] + a_{i,t} S_{i,t} -  (\mathcal{T} [\beta_{i,t}] + b_{i,t} S_{i,t} \mathcal{R}^{\top}_{i,t})\mu_{\theta_i,t+1})^{\top} \\
                     \times \mathcal{K}_{i,t+1}^{*}(\xi^*)- (\alpha_{i,t} - \beta_{i,t}\mu_{\theta_i,t+1})^{\top}\mathcal{K}_{i,t}^{*}(\xi^*) \rVert_{\mathcal{H}}
                 \end{multlined}\\
                 &=|a_{i,t}-b_{i,t}\mathcal{R}^{\top}_{i,t}\mu_{\theta_i,t+1}| \lVert  S^{\top}_{i,t} \mathcal{K}_{i,t+1}^{*}(\xi^*) \rVert_{\mathcal{H}}\\
                 &\begin{multlined}
                    =\mathcal{V}_{t+1}[|\alpha_{i,t+1} - \beta_{i,t+1}\mu_{\theta_i,t+1}|] \\
                            \times |\mathcal{K}^{**}_{i,t}(\xi_{t}) - \mathcal{K}^{*,\top}_{i,t} (\xi_{t})\mathcal{C}_{i,t} \mathcal{K}^{*}_{i,t}(\xi_{t})|
                 \end{multlined}\\
                 &= \mathcal{V}_{t+1}[|\alpha_{i,t+1} - \beta_{i,t+1}\mu_{\theta_i,t+1}|] /b_{i,t} \\ 
                 &= \mathcal{V}_{t+1}[|\alpha_{i,t+1} - \beta_{i,t+1}\mu_{\theta_i,t+1}|]/\mathcal{V}_{t+1,t+1}[\mathcal{C}_{i,t+1}]. 
            \end{aligned}
            $}
        \end{multlined}
    \end{equation}
    Following the rule for changing a data point’s position, any data point in the dataset can be moved to the last position, allowing the computations in (\ref{eq:elimination_score_proof}) to be applied to all data points, as described in (\ref{eq:elimination_score}). 
    To derive the deletion update in (\ref{eq:deletion_update}), we again start by considering the removal of the last data point (i.e., $m = t+1$) appended to the dataset to reconstruct $\alpha_{i,t}$, $\beta_{i,t}$, and $\mathcal{C}_{i,t}$. Referring to the result for $\mathcal{U}^{-}_{t+1,t+1} [\mathcal{C}_{i,t+1}]$ from (\ref{eq:K_inv}), we obtain:
    \begin{align*}
        \mathcal{U}^{-}_{t+1,t+1} [\mathcal{C}_{i,t+1}] &= \mathcal{K}^{-1}_{i,t} + \frac{\mathcal{K}_{i,t}^{-1}\mathcal{K}^{*}_{i,t} (\xi_{t})\mathcal{K}^{*,\top}_{i,t} (\xi_{t})\mathcal{K}_{i,t}^{-1}}{\mathcal{K}^{**}_{i,t}(\xi_{t}) - \mathcal{K}^{*,\top}_{i,t} (\xi_{t}) \mathcal{K}_{i,t}^{-1}\mathcal{K}^{*}_{i,t}(\xi_{t})} \\
        &= \mathcal{C}_{i,t} + \frac{\frac{\mathcal{T}^{-}_{t+1} [\mathcal{V}_{t+1}[\mathcal{C}_{i,t+1}]]^{\top}}{\mathcal{V}_{t+1,t+1}[\mathcal{C}_{i,t+1}]} \frac{\mathcal{T}^{-}_{t+1} [\mathcal{V}_{t+1}[\mathcal{C}_{i,t+1}]]}{\mathcal{V}_{t+1,t+1}[\mathcal{C}_{i,t+1}]}}{1/\mathcal{V}_{t+1,t+1}[\mathcal{C}_{i,t+1}]}.
    \end{align*}
    From this, we can directly obtain the result in (\ref{eq:deletion_update_C}) to compute  $\mathcal{C}_{i,t}$. Also, based on the definitions in (\ref{eq:GP_recursive_update_alpha}) and (\ref{eq:GP_recursive_update_beta}), one can obtain
    \begin{align*}
        \mathcal{T}^{-}_{t+1} [\alpha_{i,t+1}] &= \alpha_{i,t} + \mathcal{T}^{-}_{t+1} [a_{i,t} S_{i,t}] = \alpha_{i,t} + a_{i,t} \mathcal{K}_{i,t}^{-1}\mathcal{K}^{*}_{i,t} (\xi_{t}) \\
        &= \alpha_{i,t} + \mathcal{V}_{t+1}[\alpha_{i,t+1}] \frac{\mathcal{T}^{-}_{t+1} [\mathcal{V}_{t+1}[\mathcal{C}_{i,t+1}]]^{\top}}{\mathcal{V}_{t+1,t+1}[\mathcal{C}_{i,t+1}]},
    \end{align*}
    \begin{align*}
        \mathcal{T}^{-}_{t+1} [\beta_{i,t+1}] &= \beta_{i,t} + \mathcal{T}^{-}_{t+1} [b_{i,t} S_{i,t} P^{\top}_{i,t}] \\
        &= \beta_{i,t} + \mathcal{K}_{i,t}^{-1}\mathcal{K}^{*}_{i,t} (\xi_{t}) b_{i,t} P^{\top}_{i,t} \\
        &= \beta_{i,t} + \frac{\mathcal{T}^{-}_{t+1} [\mathcal{V}_{t+1}[\mathcal{C}_{i,t+1}]]^{\top}}{\mathcal{V}_{t+1,t+1}[\mathcal{C}_{i,t+1}]}  \mathcal{V}_{t+1}[\beta_{i,t+1}].
    \end{align*}
    To catch $\alpha_{i,t}$ and $\beta_{i,t}$ from the above expressions, we obtain the update equations in (\ref{eq:deletion_update_alpha}) and (\ref{eq:deletion_update_beta}). Following the rule for changing a data point’s position, any point in the dataset can be shifted to the final position, enabling the above computations to be carried out for all data points, as presented in (\ref{eq:deletion_update_alpha}) and (\ref{eq:deletion_update_beta}).

\section*{Acknowledgment}
This work was supported by the Research Foundation Flanders (FWO) under SBO grant no. S007723N.

\ifCLASSOPTIONcaptionsoff
  \newpage
\fi
\small{
\bibliographystyle{IEEEtran}
\bibliography{Refs}
}

\end{document}